\documentclass[10pt,twoside,final]{siamltex}
\usepackage{amsfonts,epsfig,amsmath,bm,enumerate,}
\newtheorem{remark}[theorem]{Remark}

\newtheorem{example}[theorem]{Example}
\newtheorem{examples}[theorem]{Examples}
\newtheorem{construction}[theorem]{Construction}

\def\brho{\mbox{\boldmath$\rho$}}
\def\ba{\mbox{\boldmath$a$}}
\def\bb{\mbox{\boldmath$b$}}
\def\bc{\mbox{\boldmath$c$}}
\def\bd{\mbox{\boldmath$d$}}
\def\be{\mbox{\boldmath$e$}}

\def\cos{\text{\rm cos}}

\def\Re{\text{\rm Re}}

\def\Diag{\text{\rm Diag}}
\def\lim{\text{\rm lim}}

\def\tr{\text{\rm tr}}
\def\Tr{\text{\rm Tr}}

\begin{document}
\bibliographystyle{plain}
\title{ROLE OF PARTIAL TRANSPOSE AND GENERALIZED CHOI MAPS IN QUANTUM DYNAMICAL SEMIGROUPS INVOLVING SEPARABLE AND ENTANGLED STATES\thanks{
Dedicated to Professor~R.B. Bapat on the occasion of his 60th birthday.}}

\author{Ajit Iqbal Singh\thanks{INSA Honorary Scientist, The Indian National Science Academy, New Delhi 110002, India\newline\mbox{} (ajitis@gmail.com)}}
\pagestyle{myheadings}
\markboth{A. I. Singh}{Role of Partial Transpose and Generalized Choi Maps in Quantum Dynamical Semigroups}
\maketitle

\begin{abstract}
Power symmetric matrices defined and studied by R. Sinkhorn (1981)
and their generalization by R.B. Bapat, S.K. Jain and K.~Manjunatha Prasad (1999) have been utilized to give
positive block matrices with trace one possessing positive
partial transpose, the so-called PPT states.
Another method to construct such PPT states is given, it uses the form of a matrix unitarily equivalent
to its transpose obtained by S.R.~Garcia and J.E.~Tener (2012).
Evolvement or suppression of separability or entanglement of various levels for a quantum dynamical semigroup of
completely positive maps has been studied using Choi-Jamiolkowsky matrix of such maps and the famous
 Horodecki's criteria (1996).\ A Trichotomy Theorem has been proved, and examples have been given that depend mainly on
generalized Choi maps and clearly distinguish the levels of entanglement breaking.
\end{abstract}


\section{Introduction}\label{sec-1}
Quantum inseparability or entanglement plays a significant role in
quantum communication.\ The concept goes back to A.~Einstein, E.~Schr\"odinger
and their contemporaries way back in the 1930s.\ Important
practical applications have been envisaged in recent years by
computer scientists, mathematicians and physicists. Various necessary
and sufficient conditions were given by M.~Horodecki, P.~Horodecki
and R.~Horodecki \cite{hhh}, A.~Peres \cite{ap}, for instance.\ B.M.~Terhal
and P.~Horodecki \cite{th} came up
with different levels in terms of Schmidt numbers.\ E.~St{\o}rmer
\cite{es2} has strengthened and formulated the theory in the context of
operator algebras.\
The dynamics of entanglement in continuous variable open systems with
particular emphasis on Gaussian states
has also been well studied, but we will not go into that in this paper.\
We confine our attention mainly to bipartite finite-dimensional setup in this article.

The next section is devoted to the basics of separable and entangled
states and maps as well as of quantum dynamical semigroups
in forms suitable for our purpose from standard well known books, monographs, notes,
survey articles and research papers. Particular emphasis is on
the so-called Choi maps given by M.D.~Choi
(\cite{mdc}, \cite{mdc1}) and their generalisations introduced and
studied mainly by S.J.~Cho, S.-H. Kye  and S.G.~Lee in \cite{ckl}.\
These so-called generalized Choi maps have received a lot of attention and have been
developed further by many authors.\ This section also
includes a few simple new results.\ In the third section we
give methods to construct PPT states that use power symmetric matrices due to
R.~Sinkhorn \cite{sin} and their generalization by R.B.~Bapat, S.K.~Jain and K.~Manjunatha Prasad \cite{bab}.\
Techniques developed by D.~Choudhury \cite{ch} and D.~Guillot, A.~Khare and B.~Rajaratnam
\cite{gui} help to some extent.\
The fourth section is devoted to the concept of unitary equivalence of a matrix to its transpose and its general form
obtained by S.R. Garcia and J.E.~Tener \cite{gt}.\
We formulate a variant to be exploited to give more methods to construct PPT states.
Finally, in the last section, we come to our main objective of presenting a study of
evolvement or suppression of separability or entanglement
of various levels for a quantum dynamical semigroup of completely positive maps,
in particular, the levels of entanglement breaking for
a semigroup of quantum channels.

An appendix containing a poetic felicitation to R.B. Bapat and abstract of the actual expository
talk at the conference ICLAA 2014 is given at the end.

\section{Basics of separable and entangled states and maps and
quantum dynamical semigroups}

This section is divided in eight parts A to H. Simple new results and concepts can be found in parts F, G and H.
For details an interested reader can refer to standard sources such as
those specified from place to place. On
the other hand, these can simply be skipped by workers in the area.

\renewcommand{\thesubsection}{\arabic{section}.\Alph{subsection}}
\subsection{Notation}

Let $\mathbb{N}$ be the set of
natural numbers $1,2,\ldots;\mathbb{Z}$ that of integers and
$\mathbb{Z}_{+} = \mathbb{N} \cup \{0\}.$ Let $\mathbb{R}$ be the
field of real numbers and $\mathbb{C}$ that of complex numbers. For
$n \in \mathbb{N}$, let $M_n$ denote the $C^{\ast}$-algebra of $n
\times n$ complex matrices. For $1 \leq j$, $k \leq n$, let $E_{jk}$
be the elementary $n \times n$ matrix with $1$ at the $(j,k)$-th place
and zero elsewhere. Let $I_n$ denote the identity matrix present in
$M_n.$ For $A \in M_n$, $A^t$ and $A^{\ast}$(or $A^{\dagger}$) denote
the transpose of $A$ and the adjoint of $A$ respectively. Let $\tau$
denote the transpose map on $M_n$ to itself taking $A$ to $A^t.$ Let
$M_n^{+}$ denote the positive cone of positive matrices in $M_n$, viz., the set of
positive semi-definite matrices in $M_n$.\ $A$ density matrix is an $A$ in
$M_n^{+}$ with $\mbox{tr}\,A=1$, where $\mbox{tr}$ denotes the trace.
A density $\rho$ gives rise to a positive functional $\omega_{\rho}$
on $M_n$ with $\omega_{\rho} (I_n) = 1$, also called a state, given by
$\omega_{\rho}(X) = \mbox{tr}\,(\rho X)$ for $X$ in $M_n.$ In fact,
the correspondence $\rho \rightarrow \omega_{\rho}$ is bijective and
we often use the name state for $\rho$ as well.
If the state $\rho$ has rank one, then it is called a pure state, otherwise it is called a mixed state. A
pure state $\rho$ can be given by any unit vector in its range; therefore, a unit vector in $\mathbb{C}^n$
is also called a state. The context makes it clear as to which
interpretation of the term state is being meant.

\subsection{Separable and entangled states and partial transpose}
For $n, m \in \mathbb{N}$, we
consider the tensor product $H = \mathbb{C}^m \otimes \mathbb{C}^n.$
A state on $H$ can be viewed as an $m \times n$ state and can be
represented as $\rho = [A_{jk}]_{1 \leq j, k \leq m}$ with $n \times
n$ matrices $A_{jk}$ acting on $\mathbb{C}^n.$ The state $\rho$ may
be called a bipartite state.
\begin{enumerate}[(i)]
\item The density $\rho$ is said to be {\it separable} if it is in
$M_m^{+} \otimes M_n^{+}$ in the sense that $\rho = \sum\limits_{i=1}^k p_i
\rho_i \otimes \widetilde{\rho}_i$ where $\rho_i$ and
$\widetilde{\rho}_i$ are states on $\mathcal{H}_1$ and $\mathcal{H}_2$
respectively and $p_i>0$.
\begin{enumerate}[(a)]
\item Choi \cite{mdc} gave examples of non-separable states and also necessary
conditions for $\rho$ to be separable; for example, its \emph{partial transpose}
$\rho^{PT} = [A_{jk}^t]$ is a state.\ In the literature the condition is
known as the {\it Peres test} because of significant work by Peres
\cite{ap} or {\it Positive Partial transpose}
(PPT) and, thus, states that satisfy it can be called Peres states
or {\it PPTS}.

\item A pure state $\rho$ is separable if and only if each vector $\xi$
in its range is a product vector, i.e. $\xi$ has the form $\eta\otimes\zeta$
with $\eta\in \mathbb{C}^m$ and $\zeta\in\mathbb{C}^n$; in short, $\rho$ is
a pure product state. Further, a separable mixed state is a
convex combination of such pure product states.
\end{enumerate}
\item Non-separable states are called {\it entangled states.} The acronym
{\it PPTES} is also used for an entangled PPT state in the literature
with interesting applications to Quantum Information theory appearing
mainly in Physics Journals, to name a few,
\cite{cd1}, \cite{cd2}, \cite{ch2}, \cite{ch7},  \cite{lc},
\cite{ha}, \cite{ha2}, \cite{hhh}, \cite{HSR} \cite{ph}, \cite{kk},
\cite{kk3}, \cite{pr}, \cite{sk}.
\end{enumerate}

For an expository account of this and further developments one may
consult (\cite{dav}, \cite{hia}, \cite{holevo}, \cite{hol}, \cite{ha5}, \cite{gj},
\cite{no}, \cite{nc}, \cite{par}, \cite{krp2}, \cite{Sko}, \cite{wild}).

\subsection{Separability \`a la Horodecki et al}
M. Horodecki, P. Horodecki and R. Horodecki \cite{hhh} provided
necessary
and sufficient conditions for the separability of mixed states and
gave examples to illustrate them.
\begin{enumerate}[(i)]
\item {\it Theorem} (\cite{hhh}, Theorem 2) : Let
$\mathcal{H}_1,\mathcal{H}_2$ be
Hilbert spaces of finite dimension and $\rho$ a state acting on
$\mathcal{H} = \mathcal{H}_1 \otimes \mathcal{H}_2$ i.e., $\rho$ is
a linear operator acting on $\mathcal{H}$ with $\mbox{tr}\, \rho =1$
and
$\mbox{tr}\,\rho P \ge 0$ for any projection $P.$ Let
$\mathcal{A}_1$ and $\mathcal{A}_2$ denote the sets of linear
operators acting on $\mathcal{H}_1$ and $\mathcal{H}_2$ respectively.
Then $\rho$ is {\it separable} if and only if for any positive (linear) map $\Lambda :
\mathcal{A}_2 \rightarrow \mathcal{A}_1$ the operator $(I \otimes
\Lambda) \rho$ is positive.

\item(\cite{hhh}, Remark on p.5) says that one can put
$\widetilde{\Lambda}
\otimes \Lambda$ or $\widetilde{\Lambda} \otimes I$ instead of $I
\otimes \Lambda$
(involving any positive (linear maps) $\widetilde{\Lambda}: \mathcal{A}_1
\rightarrow \mathcal{A}_2$, $\Lambda: \mathcal{A}_2 \rightarrow
\mathcal{A}_1$). The same applies to the PPT condition.

\item Next, (\cite{hhh}, Theorem 3) can be reworded as: a state
acting on
$\mathbb{C}^2 \otimes \mathbb{C}^2$ or $\mathbb{C}^2 \otimes
\mathbb{C}^3$ is separable if and only if it is PPT. The proof uses
results of St{\o}rmer \cite{es} and S.~Woronowicz \cite{slw}.
\end{enumerate}

\subsection{Entanglement breaking channels, Separable and entangled
maps}

We refer mainly to M. Horodecki, P. W. Shor and Ruskai\cite{HSR},
St{\o}rmer \cite{es2} for this subsection. The theory of completely positive maps is now a folklore. One may
consult any standard book containing the topic; we mention some
sources referred to here (\cite{as}, \cite{wba}, \cite{rbh}, \cite{mdc}, \cite{man},
\cite{mdc1}, \cite{mdc2}, \cite{dav}, \cite{ek}, \cite{krp}, \cite{es}, \cite{es1}, \cite{slw}).
Ruskai, Szarek and Werner \cite{rsw} give an interesting
analysis in the simplest set-up of $2 \times 2$ matrices with
applications to Quantum Information theory using Pauli matrices.
 However, fundamentals are given in (ii) below for the sake of convenience of the reader.
\begin{enumerate}[(i)]
\item Horodecki, Shor and Ruskai \cite{HSR} study entanglement
breaking channels. A {\it quantum channel} is a stochastic map, i.e., a map on
$M_n$ to itself which is both completely positive and trace
preserving.
\begin{enumerate}[(a)]
\item A. S. Holevo \cite{holevo} introduced channels of the form $\varphi
(\rho) = \sum\limits_{k} R_k \,\Tr\, (F_k \rho)$, where each $R_k$ is
a density matrix and $\{F_k\}$ form a positive operator valued
measure POVM. The expression for $\varphi$ is called the {\it Holevo
form} in \cite{HSR}.

\item (\cite{holevo}, \cite{HSR}, Definition 1) A stochastic map
$\varphi$ is called
{\it entanglement breaking} if $(Id \otimes \varphi)$ $A$ is
separable for any density matrix $A$, i.e., any entangled density
matrix $A$ is mapped to a separable one.

\item For $m, n \in \mathbb{N}$ and a linear map $\varphi$ on
$M_n$ to $M_m$, the {\it Choi matrix} $C_{\varphi}$ for $\varphi$ is
$C_{\varphi} = \sum\limits_{j,k} \,E_{jk} \otimes \varphi
(E_{jk}) = (I d \otimes \varphi) \rho \in M_n \otimes
M_m$, where $\frac{1}{n} \rho$ is the
1-dimensional projection $\frac{1}{n}\sum\limits_{j,k} E_{jk} \otimes
E_{jk}$, the so-called basic maximally entangled state. M.D.~Choi
(\cite{man}) proved that
$\varphi$ is completely positive if and only if $C_{\varphi}$ is
positive. Physicists usually use $\frac{1}{n} C_{\varphi}$ for a
trace-preserving completely positive map $\varphi$ and call it a
Jamiolkowski state, (See for instance, \cite{rw}, \cite{wc})
following A. Jamiolkowski \cite{aj}.

\item A part of (\cite{HSR}, Theorem 4) says that the following are
equivalent for a channel $\varphi.$
\begin{enumerate}[$(\alpha)$]
\item[$(\alpha)$] $\varphi$ is entanglement breaking.
\item[$(\beta)$] $\varphi$ has the Holevo form with $F_k$ positive semi-definite.
\item[$(\gamma)$] $\frac{1}{n} C_{\varphi}$ is separable.
\item[$(\delta)$] $\psi \circ \varphi$ is completely positive for all
 positivity preserving maps $\psi.$
\item[$(\in)$] $\varphi \circ \Lambda$ is completely positive for all
 positivity preserving maps $\Lambda.$
\end{enumerate}
\end{enumerate}

\item (\cite{es2},\S 1) Let $\mathcal{A}$ be an operator system, i.e.,
a norm-closed self-adjoint linear space of bounded
operators on a
Hilbert space $\mathcal{K}$ containing the identity. Let
$\mathcal{H}$ be a Hilbert space and $\mathcal{B}(\mathcal{H})$, its
operator algebra and $\mathcal{T}(\mathcal{H})$, the space of the
trace class operators on $\mathcal{H}.$ Let $\tau$ be the transpose
map
on $\mathcal{B}(\mathcal{H})$ (respectively
$\mathcal{B}(\mathcal{K})$) with respect to some orthonormal basis for
$\mathcal{H}$ (respectively $\mathcal{K}$). At times for $a \in
\mathcal{B}(\mathcal{H})$ or $\mathcal{B}(\mathcal{K})$, $\tau (a)$
will be denoted by $a^t.$

The $BW$-topology on the space of
bounded linear maps on $\mathcal{A}$ to $\mathcal{B}(\mathcal{H})$ is
the topology of bounded pointwise weak convergence, i.e., a net
$(\varphi_{\nu})$ converges to $\varphi$ if it is uniformly bounded,
and $\varphi_{\nu} (a) \rightarrow \varphi(a)$ weakly for all $a
\in \mathcal{A}.$ Let $S(\mathcal{H})$ be the $BW$-closed cone generated
by maps of the form
$$x \rightarrow \sum_{j=1}^n \omega_j (x) a_j$$
where $\omega_j$ is a normal state on $\mathcal{B}(\mathcal{H})$ and
$a_j \in$ the positive cone $\mathcal{B}(\mathcal{H})^{+}.$ Here,
$\mathcal{A} = \mathcal{B}(\mathcal{H}).$

For the sake of convenience we recall some well-known conditions for a linear map
$\varphi$
on $\mathcal{A}$ to $\mathcal{B} = \mathcal{B}(\mathcal{H}).$ Here $r
\in
\mathbb{N}$, and $\varphi$ is a {\it $\ast$-map} in the sense that
$\varphi(x^{\ast}) = \varphi (x)^{\ast}$ for $x \in \mathcal{A}.$
\begin{enumerate}[(a)]
\item The map $\varphi$ is said to be {\it $r$-positive} if the map
 $\varphi_r = \varphi \otimes\, {\rm Id}: \mathcal{A} \otimes M_r \rightarrow
\mathcal{B} \otimes M_r$ is positive.

\item The map $\varphi$ is said to be {\it completely positive} if
$\varphi$ is $r$-positive for all $r \in \mathbb{N}.$
\item The map $\varphi$ is said to be {\it $r$-copositive (respectively,
completely co-positive)}, if $\tau \circ \varphi$ is $r$-positive
(respectively, completely positive).

\item The map $\varphi$ is said to be a {\it Schwarz} map if
$$\varphi(x^{\ast} x) \ge \varphi (x^{\ast}) \varphi (x)
\,\,\mbox{for}\,\, x \in \mathcal{A} \,\,\mbox{with}\,\, x^{\ast} x
\in \mathcal{A}\,.$$
\item The map $\varphi$ is said to be $r$-{\it Schwarz} if $\varphi_r$
is a
Schwarz map i.e.,
\begin{align*}
\varphi_r \left (\left [x_{jk} \right ]^{\ast} \left [x_{jk} \right ]
\right)
&\ge \varphi_r \left (\left [x_{jk} \right ]^{\ast} \right) \varphi_r \left (\left [x_{jk} \right ] \right)\\
&\qquad \mbox{for}\,\, x_{jk} \in \mathcal{A}, 1 \leq j, k \leq r \,\,\mbox{with} \,\, x_{jk}^{\ast} x_{pq} \in \mathcal{A}\\
&\qquad \mbox{for} \,\, 1 \leq j,k,p,q \leq r.
\end{align*}
\item For a $C^{\ast}$-algebra $\mathcal{A}$, and $\varphi$ unital in
the sense that $\varphi$ takes the identity of $\mathcal{A}$ to that of
$\mathcal{B}(\mathcal{H})$, $\varphi$ is 2-positive
implies that $\varphi$ is a Schwarz map. As a consequence, a unital
$\varphi$
is completely positive if and only if $\varphi$ is $r$-Schwarz for
each $r.$

\item For a unital $\varphi$, the inequality in (d) above is satisfied
for normal elements $x$ when $\varphi$ is positive.

\item Many more assorted inequalities hold for positive maps (see
\cite{mdc1}, for instance).

\item For $\varphi$ satisfying any condition as in (a) to (e), $\psi$
completely positive on $\mathcal{B}(\mathcal{H})$ to itself, and
$\Lambda$ completely positive on $\mathcal{A}$ to itself, the maps
$\varphi \circ \Lambda$ and $\psi \circ \varphi$ satisfy the
corresponding condition.

\item Similar terminology as in (c) above applies to other conditions
like (d) and~(e).
\end{enumerate}

\item (\cite{es2}, Lemma 1) sets up an isometric isomorphism $\varphi
\rightarrow \widetilde{\varphi}$ between the set
$\mathcal{B}(\mathcal{A},\mathcal{B}(\mathcal{H}))$ of bounded linear
maps of $\mathcal{A}$ into $\mathcal{B}(\mathcal{H})$ and the dual
$(\mathcal{A} \widehat{\otimes}\mathcal{T}(\mathcal{H}))^{\ast}$ of
the projective tensor product of $\mathcal{A}$ and
$\mathcal{T}(\mathcal{H})$ given by
$$\widetilde{\varphi} (a \otimes b) = \mbox{Tr}\, (\varphi (a) b^t) $$
where $\mbox{Tr}$ denotes the usual trace on $\mathcal{B}(\mathcal{H})$ taking
the value 1 on minimal projections. Furthermore, $\varphi$ is a positive linear
operator if and only if $\widetilde
{\varphi}$ is positive on the cone $\mathcal{A}^{+}
\widehat{\otimes} \mathcal{T}(H)^{+}$ generated by operators of the
form $a \otimes b$ with $a$ and $b$ positive.

\item As noted in \cite{es2}, p.~2305, it follows from (\cite{es1},
Theorem 3.2) that $\varphi$ is completely positive if and only if
$\widetilde{\varphi}$ is positive on the cone $(\mathcal{A}
\widehat{\otimes} \mathcal{T}(\mathcal{H}))^{+}$, the closure of
the positive operators in the algebraic tensor product $\mathcal{A}
\odot \mathcal{T}(\mathcal{H}).$

\item A positive linear functional $\rho$ on $\mathcal{A}
\widehat{\otimes}
\mathcal{T}(\mathcal{H})$ is said to be {\it separable} if it belongs
to the
norm closure of positive sums of states of the form $\sigma \otimes \omega$
where $\sigma$ is a state of $\mathcal{A}$ and $\omega$ a normal
state of
$\mathcal{B}(\mathcal{H}).$ Otherwise $\rho$ is called {\it
entangled}.

\item A part of
(\cite{es2}, Theorem 2) says that the following are equivalent for a
$\varphi \in \mathcal{B}(\mathcal{A}, \mathcal{B} (\mathcal{H})).$
\begin{enumerate}
 \item[(a)] $\widetilde{\varphi}$ is a separable positive linear functional.

 \item[(b)] $\varphi$ is a BW-limit of maps of the form $x \rightarrow
\sum\limits_{j=1}^{n} \omega_j (x) b_j$ with $\omega_j$ a state of
$\mathcal{A}$
and
$b_j \in \mathcal{B}(\mathcal{H})^{+}.$
\end{enumerate}

\item \textsc{Definition.} A completely positive map $\varphi$ in
$\mathcal{B}(\mathcal{A},
\mathcal{B}(\mathcal{H}))$ will be called {\it separable} (respectively, {\it entangled)} if $\widetilde{\varphi}$ is so. A
separable map may also be called \emph{entanglement breaking},
in view of (i)(d) above for special $\varphi$'s, viz., channels, and
we do so at times in what follows.

\item (\cite{es2}, Corollary 3) can now be reworded as : Let
$\mathcal{H}$ be separable
and $\varphi \in \mathcal{B}(\mathcal{A}, \,
\mathcal{B}(\mathcal{H}))$ be
positive. If $\varphi
(\mathcal{A})$ is contained in an abelian $C^{\ast}$-algebra then
$\varphi$ is
separable.

\item We shall say a positive linear functional $\rho$ on $\mathcal{A}
\widehat{\otimes}
\mathcal{T}(\mathcal{H})$ is {\it PPT} (i.e. satisfies the Peres
condition)
if
$\rho \circ (Id \otimes \tau))$ is positive. In line with (ii)
above $\varphi
\in \mathcal{B}(\mathcal{A}, \mathcal{B}(\mathcal{H}))$ will be said
to be {\it PPT}
if $\widetilde{\varphi}$ is so. (\cite{es2}, Proposition 4) can now be
interpreted
as: $\varphi$
is PPT if and only if $\varphi$ is both completely positive and completely
co-positive.

\end{enumerate}

We are now ready to prove a simple result in line with item (i)(d)
above.

\begin{theorem}
Let $\varphi \in \mathcal{B}(\mathcal{A},
\mathcal{B}(\mathcal{H}))$ be a separable map. Then for any positive
bounded map $\psi$ on $\mathcal{B}(\mathcal{H})$ to itself and any positive unital map
$\Lambda$ on $\mathcal{A}$ to itself, the maps $\psi \circ \varphi$
and $\varphi \circ \Lambda$ are both separable.
\end{theorem}

\proof We note that for a state $\omega$ of $\mathcal{A}$
and $b \in \mathcal{B}(\mathcal{H})^{+}$, $\omega \circ \Lambda$ is
a
state of
$\mathcal{A}$ and $\psi(b) \in \mathcal{B}(\mathcal{H})^{+}.$ We can now apply
the item 2.D(vi) above for $\varphi$ and conclude that maps $\varphi
\circ \Lambda$
and $\psi \circ \varphi$ satisfy the condition (b) in 2.D(vi). By
2.D(vi) we obtain
that $\varphi \circ \Lambda$ and $\psi \circ \varphi$ are
separable.\endproof

\begin{corollary} Let $\varphi$ be a completely positive map on
$\mathcal{A}$ to $\mathcal{B}(\mathcal{H}).$ If $\varphi$ is separable then for
each positive bounded map $\psi$ on $\mathcal{B}(\mathcal{H})$ to
itself and each positive
unital map $\Lambda$ on $\mathcal{A}$ to
itself, $\psi \circ \varphi$ and
$\varphi \circ \Lambda$ are
completely positive.
\end{corollary}

\subsection{Horodecki's-St{\o}rmer Theorem}
Let $m, n \in \mathbb{N}$, $\mathcal{A} = M_n$, $\mathcal{H} = \mathcal{C}^m$,
$\varphi \in \mathcal{B}(A, \mathcal{B}(\mathcal{H}))$, and
$C_{\varphi}$ the Choi matrix for
$\varphi.$ The map $\varphi^t = \tau\circ \varphi\circ \tau$
(where $\tau$ is
the transpose map in either $M_n$ or $M_m$) is completely positive if and only
if $\varphi$ is so. (\cite{es2}, Lemma
5) says
that $C_{\varphi^{t}}$ is the density matrix for
$\widetilde{\varphi}.$

\medskip
St{\o}rmer continues with his study and gives, amongst other things,
his infinite-dimensional extension of Horodecki's Theorem, which we
may call Horodecki's-\break \mbox{St{\o}rmer~Theorem}, and also methods to construct
PPTES. We shall not go into details here (cf.\ \cite{ha4}, \cite{joh3}, \cite{yl}, \cite{sko1}).

\subsection{Pure Product states and Schmidt number}
\begin{enumerate}[(i)]
\item P.~Horodecki \cite{ph} proved that a separable state on the Hilbert space
$\mathcal{H}=\mathcal{H}_1 \otimes \mathcal{H}_2$ (with $\mathcal{H}_1$ and
$\mathcal{H}_2$ both finite-dimensional) can be written as a convex combination
 of $N$ pure product states with $N \leq (\mbox{dim} \mathcal{H})^2$ and gave a
new separability criterion in terms of the range of the density matrix. This
was carried further in different ways by several authors (cf.~\cite{lc},
\cite{ha}, \cite{vsw}, \cite{kk}, \cite{krp1}, \cite{th}, \cite{w}).

\item K. R. Parthasarathy \cite{krp1} called a subspace of $\mathcal{H}
=
\mathcal{H}_1 \otimes\mathcal{H}_2 \otimes \cdots \otimes\mathcal{H}_r$
\emph{completely entangled} if it contains no non-zero
product vector of the form $u_1 \otimes u_2
\otimes \cdots \otimes u_r$ and gave a concrete example of a space attaining the maximal
dimension of such a space. He also introduced a more delicate notion
of perfectly entangled subspace
for a multipartite quantum system. The notion of completely entangled
subspaces is related to notions of unextendible product bases and
uncompletable product bases, which are well studied by C.H. Bennett et al \cite{benn},
D.P. Di Vincenzo et al
\cite{vmsst} and further studied by N. Alon and L. Lov\'asz \cite{al},
A.O. Pittenger \cite{aop}, B.V.R. Bhat \cite{bvrb}, L. Skowronek
\cite{sk}, N. Johnston \cite{joh1,joh2}, R. Sengupta, Arvind and the author \cite{sen} and many more authors.
They have been used to construct entangled PPT
densities on one hand and non-PPT ones on the other.

\item B.M.~Terhal and P.~Horodecki \cite{th} extended the notion of
the Schmidt rank or
number of a pure state to the domain of bipartite density matrices.
\begin{enumerate}[(a)]
\item The Schmidt rank (or number) of a pure bipartite state $\xi$
is the number of non zero co-efficients
 in the essentially unique Schmidt form
 $\xi=\sum_j \lambda_j \eta_j\otimes \zeta_j$ with $\eta_j$'s and $\zeta_j$'s forming
 orthonormal sets in their respective spaces.
 For a mixed bipartite state $\rho$,
 its Schmidt rank is the maximum Schmidt rank (or number)
 in an optimal pure state decomposition of $\rho$.

\item To motivate
our next notions, we quote their characterization viz., Theorem~1 \cite{th}.
\end{enumerate}

Let $\rho$ be a density matrix on $\mathcal{H}_n \otimes \mathcal{H}_n$, i.e., $\mathbb{C}^n\otimes \mathbb{C}^n$.\ The
density matrix has the Schmidt number at least $r + 1$ if and only if there
exists an $r$-positive linear map $\Lambda: M_n \rightarrow M_n$ such that $(I
\otimes \Lambda) (\rho) \not\ge 0.$
\end{enumerate}

\newpage
\begin{definition}\rm
Let $\mathcal{A}$, $\mathcal{B}(\mathcal{H})$
etc. be as in 2.D and $\varphi \in\mathcal{B}
(\mathcal{A},\mathcal{B}(\mathcal{H}))$ be a completely positive map. We say it
has \emph{Schmidt number at least $r + 1$} if either there exists an
$r$-positive linear map $\psi :\mathcal{B} (\mathcal{H}) \rightarrow
\mathcal{B}(\mathcal{H})$ such that $\psi \circ
\varphi$ is not completely positive or there exists an $r$-positive linear map
$\Lambda : \mathcal{A} \rightarrow \mathcal{A}$ such that $\varphi
\circ \Lambda$ is not completely positive.

Remarkable progress on the topic has been made in recent years as can be gauged
from (\cite{ch1}, \cite{ch3}, \cite{ch4}, \cite{ch5}, \cite{ch6}, \cite{ch7},
\cite{ch8}, \cite{ha5}, \cite{joh3}, \cite{kk1}, \cite{kk2}, \cite{yl},
\cite{mm19}, \cite{mm20}, \cite{ar}, \cite{sko2})
The name ``partially entanglement breaking '' has been associated with
maps having Schmidt number $<n$ by some authors.\
We shall not go into details in this paper.
\end{definition}

\subsection{Quantum dynamical semigroups}
We may refer to any
standard source for this folklore material, particularly (\cite{ce},
\cite{en}, \cite{ek},
\cite{vg}, \cite{gl}, \cite{krp}, \cite{wc}) mentioned in the list of
references, if we like, rather than original sources.
\begin{enumerate}[(i)]
\item Algebraically speaking, a dynamical system is a family
$(T(t))_{t \ge 0}$ or, for short $(T_t)_{t \ge 0}$, of mappings on
a set $\mathcal{X}$ satisfying
\begin{align*}
T (t+s) &= T(t) T(s) \,\,\mbox{for all} \,\, t, s \ge 0 \\
T(0) &= Id.
\end{align*}
In fact even if we just confine our attention to the first condition,
we have
$T_0^2 = T_{0}$, the range $\mathcal{R}_0$ of $T_0$ contains the range
$\mathcal{R}_{t}$
of $T_t$ for all $t$, and $T_0$ restricted to $\mathcal{R}_{0}$ is
the identity. So we may replace $\mathcal{X}$ by $\mathcal{R}_{0}$,
and then the second condition holds for $(S_t)_{t\ge 0}$, where $S_t=T_t |
\mathcal{R}_0.$ Then $T_t = S_t T_0 = T_t T_0 = T_0 T_t$ and $S_t
S_s = S_s S_t = S_{t+s}.$ Thus, $(S_t)_{t \ge 0}$ is a dynamical
system and we call $(T_t)_{t\ge0}$ a {\it $T_0$-constricted dynamical
system}.
\item Usually $\mathcal{X}$ is taken to be a Banach space, $T_t$ a bounded linear
operator on $\mathcal{X}$ for each $t$ and the system to be strongly continuous
on $\mathcal{X}$. Then $(T_t)_{t\ge 0}$ is called a {\it strongly continuous}
({\it one parameter}) {\it semigroup} or a {\it $C_0$-semigroup}.
Again, if we relax the condition $T_0=Id$, we call $(T_t)_{t\ge0}$ a
{\it $T_0$-constricted $C_0$-semigroup.} We note that the continuity
of $T_0$ and the fact that $T_0 x = x$ for $x$ in $\mathcal{R}_0$
forces $T_0x = x$ for $x$ in the closure $\bar{\mathcal{R}}_0$ of
$\mathcal{R}_0.$ This, in turn, gives that $\bar{\mathcal{R}}_0
\subset \mathcal{R}_0.$ Hence $\mathcal{R}_0$ is closed and,
therefore, a Banach space.
\item When $\mathcal{X}$ is an operator system $\mathcal{A}$, and maps $T_t$
satisfy conditions like those in the item 2.D(ii) (a) to (e) above,
we term a $C_0$-semigroup
$(T_t)_{t\ge 0}$ as a {\it quantum dynamical semigroup or system.} In
practice, $T_t$'s are
all taken to be completely positive and $\mathcal{A}$ to be $M_n$ or
$\mathcal{B}(\mathcal{H}).$ Once again, the term $T_0$-constricted
quantum dynamical semigroup will be used when we relax the
condition $T_0=Id.$
\item For a $C_0$-semigroup $(T_t)_{t\ge 0}$, the \emph{infinitesimal
generator} $L$ is the operator which has the domain
$$D(L) = \left \{x \in \mathcal{X} :\underset{t \rightarrow
0 +}{\lim}
\frac{1}{t} (T_t x-x) \,\,\mbox{exists} \right \} $$
and is given by $Lx = \underset{t \rightarrow 0+}{\lim}
\frac{1}{t} (T_t x -x)$ for $x$ in $D(L)$. Then $L$
 is a closed and densely-defined linear operator that determines the semigroup
uniquely.
\item For an $A \in \mathcal{B}(\mathcal{X})$, $T_t = \exp (tA)$, $L$
coincides
with $A.$ For this reason $T_t$ as in (iv) above is written as
$\exp(tL) = e^{tL}$
as well.

\item
When $T_t$'s
satisfy any of the conditions in 2.D(ii) (a) to (e) or
corresponding ``co'' parts as indicated in 2.D(ii) (j) above, $L$
satisfies a
corresponding variant of the condition. Fundamental theoretical work in this
direction is by V. Gorini, A. Kossakowski and E. C. G. Sudarashan
\cite{vg} and
G. Lindblad \cite{gl}, though history can be traced back to specific
irreversible processes or quantum stochastic processes of open systems
by many
like R. V. Kadison or E. B. Davies. For further basic developments
one can see \cite{ce}, \cite{ek} and \cite{krp}.

\item
It follows from the proof of and the Proposition itself on p.73
\cite{en}
that if there exists some $t_0 > 0$ such that $T(t_0)$ is invertible,
then
\begin{enumerate}[(a)]
\item for $0 \leq t < t_0$, $T(t_0)= T(t_0-t)T(t) = T(t) T (t_0 - t)$ and
for $t = nt_0 +
s$ for $n \in \mathbb{N}$, $s \in [0, t_0)$, $T(t) = T(t_0)^n T(s)$,
and therefore, $T(t)$ is invertible for all $t \ge 0;$ $T(t)^{-1} = T(t_0-
t) T(t_0)^{-1}$ for $0 \leq t < t_0$ and $T(t)^{-1} =
T(s)^{-1}
T(t_0)^{-n}$ for $t = nt_0 + s$ for $n \in \mathbb{N}$, $s \in [0,
t_0);$

\item $(T(t))_{t \ge 0}$ can be embedded in a group $(T(t))_{t \in
\mathbb{R}}$ on $\mathcal{X}.$
\end{enumerate}
\end{enumerate}

\begin{theorem} Let $(T_t)_{t \ge 0}$ be a quantum dynamical
system of completely positive maps. If there exists some $t_0 > 0$ such that
$T(t_0)$ is invertible and $T(t_0)^{-1}$ satisfies any of the conditions 2.D(ii) (a)
to (c) then each $T(t)^{-1}$ satisfies the corresponding condition.
\end{theorem}

\proof This is obvious from (vii)(a) just above.\endproof

\subsection{Separability and entanglement of generalized Choi maps and use in Quantum Information theory}
We begin with a general set-up.

\begin{enumerate}[(i)]
\item \textbf{A foliation:}\ Let $D_n$ be the linear span of $\{
E_{jj} : 1 \leq j \leq
n\}$ and $F_n$ the linear span of $\{E_{jk} : 1 \leq j \neq k \leq n
\}.$ We note that as a linear space $M_n = D_n \oplus F_n.$ Also any
linear map $\Lambda$ on $M_n$ to itself can be expressed in the form
$\left [\begin{array}{cc} \Lambda_{11} & \Lambda_{12} \\
\Lambda_{21} & \Lambda_{22} \end{array} \right ]$ where $\Lambda_{11}
: D_n \rightarrow D_n,$ $\Lambda_{12} : F_n \rightarrow D_n,$
$\Lambda_{21} : D_n \rightarrow F_n$ and $\Lambda_{22} : F_n
\rightarrow F_n$ are linear maps.
\begin{enumerate}[(a)]
\item Let $C_{\Lambda}$ be the
so-called Choi
matrix of map $\Lambda.$ It is given by the block matrix $\left
[\Lambda (E_{jk}) \right
]$ written as an $n^2 \times n^2$ matrix with entries in
$\mathbb{C},$ in fact. The diagonal of $C_{\Lambda}$ is same as the
diagonal of the
block matrix with $\Lambda_{11} (E_{jj})$ at the $jj^{\mbox{th}}$
block. As a consequence, $\mbox{tr} (C_{\Lambda}) =
\sum\limits_{j=1}^n$ $\mbox{tr} \Lambda_{11} (E_{jj}) = \mbox{tr}\,
\Lambda_{11} (I_n) = \mbox{tr} \,\Lambda (I_n).$ So $C_{\Lambda}$ is
a density matrix if and only if $\Lambda$ is completely positive with
$\mbox{tr}\,\Lambda_{11} (I_n)=1.$ See (\cite{mdc}, \cite{man}, \cite{mdc1},
\cite{mdc2})
for more details.

\item We consider the class $\mathcal{L}$
of maps $\Lambda$ with $\Lambda_{12} = 0$ and $\Lambda_{21} = 0$ and
write $\Lambda = \Lambda_1 \oplus \Lambda_2$ with $\Lambda_1 =
\Lambda_{11}$ and $\Lambda_2 = \Lambda_{22}.$ Addition and product of
maps in $\mathcal{L}$ is component-wise and as a consequence, for
$\Lambda \in \mathcal{L},$ $e^{\Lambda} = e^{\Lambda_{1}} \oplus
e^{\Lambda_{2}}.$

\item A large number of examples in the study of
positive, $k$-positive and completely positive maps, of dynamical
semigroups and of separability, entanglement and Schmidt number for
density matrices are in $\mathcal{L}$. One may observe this tendency in
\cite{ckl}, \cite{mdc}, \cite{mdc1}, \cite{mdc2}, \cite{ch1}, \cite{ch2},
\cite{ch3}, \cite{ch4}, \cite{ch5}, \cite{ch6}, \cite{ch7}, \cite{ch8},
\cite{lc}, \cite{ha1}, \cite{ha3}, \cite{hhh}, \cite{ph}, \cite{bk},
\cite{shk}, \cite{mm19}, \cite{ho21}, \cite{pr}, \cite{rsw}, \cite{ar},
\cite{sen2}, \cite{sen1}, \cite{th}, for instance.

Such maps have been tested for more properties like separability,
entanglement and Schmidt numbers and used to construct
Entanglement witnesses, Entanglement
breaking or partially entanglement breaking channels by some of
them. We shall not go into details here.
\item
Quite often
$\Lambda_2$ is just in the
one-dimensional linear space spanned by $\mathcal{I}_{F_{n}}$, the
identity operator on $F_n$, or else in the two dimensional algebra
generated by $\mathcal{I}_{F_{n}}$ and the restriction $\tau_{F_{n}}$
of the transpose map $\tau$ on $M_n.$ Also $\Lambda_1$'s are usually
taken to be upper (or lower) triangular matrices (cf.\cite{bk}) or
matrices with rows being just permutations of each other (\cite{ckl},
\cite{mdc1}).
\end{enumerate}
We confine our attention to only one class and give below relevant known results
or minor riders for further use in Section 5.

\item \textbf{Generalized Choi maps:}\
We reformulate generalized Choi maps as presented in
Cho, Kye and Lee \cite{ckl} and also give little variants of them below.
\begin{enumerate}[(a)]
\item For $a,b,c \in \mathbb{C},$ let $$D(a,b,c) =
\begin{pmatrix} a & b & c \\ c & a & b \\ b & c & a\end{pmatrix}.$$

The set $\mathcal{D} = \{ D(a,b,c):a,b,c \in \mathbb{C}\}$ is a
commutative semigroup with identity $D(1,0,0)=I_3$ simply because $D
(a^{\prime},b^{\prime},c^{\prime})$ $D(a,b,c)= D (a^{\prime}a +
b^{\prime}c + c^{\prime} b, c^{\prime} c + a^{\prime} b + b^{\prime}
a, b^{\prime} b + c^{\prime} a +a^{\prime} c ).$ Further $\mathcal{D}
\cap G L (3, \mathbb{C})$ is a subgroup of the general linear group
$GL(3,
\mathbb{C}).$ To see
this it is enough to note that
\begin{align*}
D (a,b,c) D (a^2 -bc, c^2 -ab, b^2-ac)
&=(a^3+b^3+c^3-3abc)I_3 \\
&= \det D(a,b,c)I_3.
\end{align*}
\item For $a,b,c \in \mathbb{R}_{+} = [0, \infty),$ the map $\Phi
[a,b,c]$
defined on p.214 \cite{ckl} is the same as $D(a-1, b,c) \oplus
(-\mathcal{I}_{F_{3}}).$ We prefer to consider for $a,b,c \in \mathbb{R}_{+},$
the variants
$$\rho [a,b,c] = \Phi [a+1, b, c] = D(a,b,c) \oplus (-\mathcal{I}_{F_{3}})$$
and generalizations,
\begin{align*}
\rho [a,b,c,d] &= D (a,b,c) \oplus d \mathcal{I}_{F_{3}}, \quad \mbox{and}
\end{align*}
\begin{align*}
\tau [a,b,c,d] &= D (a,b,c) \oplus d \, \tau_{F_{3}} =\rho[a,b,c,d]\,\tau = \,\tau \rho[a,b,c,d].
\end{align*}

\item We note that $\rho[a,b,c,d]$ is unital if and only if $a+b+c=1$
if and only if $\tau[a,b,c,d]$ is unital. Also $\rho[a,b,c,d]$ is
trace-preserving if and only if $a+b+c=1$ if and only if
$\tau[a,b,c,d]$ is trace-preserving.

The map $\rho[a,b,c,d]$ is a $\ast$-map if and only if $a,b,c,d
\in \mathbb{R}$ if and only if $\tau[a,b,c,d]$ is a $\ast$-map.

Finally, if $\rho[a,b,c,d]$ or $\tau[a,b,c,d]$ is a positive map then
$a,b,c$ are all non-negative simply because the image of $E_{11}$ is
$\begin{bmatrix}
a & 0 & 0 \\ 0 & c & 0 \\ 0 & 0 & b \\
\end{bmatrix}.$

\item We note that $\rho[1, 0, \mu]$ with $\mu \ge 1$ is the same as
the map $\Phi$ in (\cite{mdc1}, Appendix B, Example).

\end{enumerate}

For the sake of convenience we recast the results in \cite{ckl} in our notation.

\item \textbf{Properties of $\rho[a,b,c]$, $\rho[a,b,c,d]$ and $\tau[a,b,c,d]$}.

Let $a,b,c \in \mathbb{R}_{+},$ $d \in \mathbb{R}.$
\begin{enumerate}[(a)]
\item By Theorem 2.1 \cite{ckl}, the map $\rho[a,b,c]$ is
positive if and only if $a+b+c \ge 2$ together with $bc \ge (1-a)^2$ in case $0 \leq a \leq 1.$

\item
By \cite{ckl}, Lemma~3.1, the map $\rho [a,b,c,d]$ is
completely positive if and only if $a \ge d$ and $a \ge -2d.$

In particular, $\rho[a,b,c]$ is completely positive if and only
if $a \ge 2.$ This is Proposition 3.2 \cite{ckl}.

\item
By \cite{ckl} Lemma 3.1, second part, $\tau [a,b,c,d]$ is
completely positive if and only if $bc \ge d^2.$ As a consequence,
$\rho [a, b, c, d]$ is positive if $b c \ge d^2.$

In particular, $\tau \circ \rho [a,b,c]$ is completely
 positive (i.e. $\rho[a,b,c]$ is completely copositive) if and only
if
$bc \ge 1.$ This is a part of Proposition 3.3 of \cite{ckl}.

\item
Theorem 3.4 \cite{ckl} gives that for $0 \leq a <
2,$ $\rho [a,b,c]$ is decomposable if and only if $bc \ge
(1-\frac{a}{2})^2.$

This combined with (b) and (c)
gives that $\rho[0,b,c]$ is not completely positive and it is
decomposable if and only if $bc \ge 1$ if and only if it is
completely copositive.

Further, for $0 < a < 2$ and $bc < 1$, $\rho[a,b,c]$ is neither completely
positive nor completely copositive, but is nevertheless decomposable
if and only if $bc \ge (1-\frac{a}{2})^2.$

\item Theorem 4.2 [11] says that $\rho[a,b,c]$ is 2-positive if and only if either $a\ge 2$ or
$1\le a <2$ and $bc=(2-a)(b+c) > 0$.
Further, taking $a=1$, $b=2=c$, we get the well-known example of a 2-positive map that is not
completely positive given by Choi in [12] as remarked on p.~214 [11] as well.

\end{enumerate}
\end{enumerate}

\begin{remark}\rm
Corollary~2.2 combined with the above remarks can give us
a multitude of completely positive maps and
states that are PPT, non-PPT or PPTE. We illustrate this by recording a
few special cases for further use.
\end{remark}

\begin{theorem} Let $a,b,c,d \in \mathbb{R}$ with
$a,b,c \ge 0.$ Let $\varphi = \rho [a,b,c,d]$ and $\psi = \tau
[a,b,c,d].$
\begin{enumerate}[\rm(i)]
\item If $d =0,$ then $\varphi = \psi$ is a separable map.
\item Let $d>0.$
\begin{enumerate}[\rm(a)]
\item If $a \ge d$ but $bc < d^2,$ then $\varphi$ is a non-PPT completely positive map.
\item If $bc \ge d^2$ but $a < d,$ then $\psi$ is a non-PPT completely positive map.
\item If $a \ge d,$ $bc \ge d^2$, then $\varphi$ and $\psi$ are PPT maps.
\item If $a \ge d,$ $bc \ge d^2$ but $a+b < 2d$ or $a+c < 2d$, then $\varphi$ and $\psi$ are PPTE maps.

In particular, this is so if $a=1=d,$ $0 < b < 1,$ $c=\frac{1}{b}.$

\item If  $a+2(b+c)<2d\le 2a$ then $\varphi$ has Schmidt number 3 and it is non-PPT. In particular, it is so if $b+c=1-a$,
and, $\frac{2}{3}<a\le 1$ together with $1-\frac{a}{2}<d\le a$, or, equivalently, $\frac{1}{2}<d\le 1$ together with $d\le a\le 1$ for $d>\frac{2}{3}$ whereas
$2(1-d)<a\le 1$ for $\frac{1}{2}<d\le \frac{2}{3}$.

\end{enumerate}
\item Let $d < 0.$
\begin{enumerate}[\rm(a)]
\item If $a \ge -2d= 2|d|$ but $bc < d^2$, then $\varphi$ is a non-PPT completely positive map.
\item If $bc \ge d^2$ but $a < 2|d|,$ then $\psi$ is a non-PPT completely positive map.
\item If $a \ge 2|d|,$ $bc \ge d^2$, then $\varphi$ and $\psi$ are PPT maps.
\end{enumerate}
\end{enumerate}
\end{theorem}

\proof Parts (i), (ii) (a), (ii)(b), (ii)(c),
(iii)(a), (iii)(b) follow immediately from (iii) above. For
(ii)(d) we consider the map $\xi= \rho[a^{\prime}, b^{\prime},
c^{\prime}] \equiv \rho[a^{\prime}, b^{\prime},
c^{\prime}, -1]$ with $a^{\prime}, b^{\prime}, c^{\prime}$ to be
suitably chosen yet to be specified. We have
$$\xi \varphi = \rho [a^{\prime}a + b^{\prime}c+c^{\prime}b,
c^{\prime}c + a^{\prime}b+b^{\prime}a, b^{\prime}b +
c^{\prime}a+a^{\prime}c, - d].$$

By item (iii)(b), $\xi \varphi$ is completely positive if and
only if $a^{\prime}a + b^{\prime}c + c^{\prime} b \ge 2d.$ If $a + b
< 2d$ then we take $a^{\prime} = 1 = c^{\prime}$ and $b^{\prime} =0.$
Then $a^{\prime}+ b^{\prime} + c^{\prime} =2,$ $b^{\prime} c^{\prime}
= 0 =(1-a^{\prime})^2$ and $a^{\prime} a + b^{\prime}c + c^{\prime} b
= a+b < 2d.$ So by Item (iii)(a) and (b),  $\xi$ is
positive and $\xi \varphi$ is not completely positive. Similarly, in
case $a + c < 2d,$ we take $a^{\prime}=1=b^{\prime}$ and
$c^{\prime}=0$ and conclude that $\xi$ is positive but $\xi \varphi$ is
not completely positive. By Corollary~2.2,
$\varphi$ and $\psi$ are not separable. So $\varphi$ and $\psi$ are
PPTE maps.

For (ii)(e), we may consider $\eta= \rho[1,2,2]$, which as noted in item (iii)(e) above,  is
2-positive but not completely positive. We have
$\eta \phi=\rho[a+2(b+c), 2c+b+2a, 2b+c+2a, -d]$, which , by item (iii)(b) can not
be completely positive, simply because $a+2(b+c) < -2(-d)$.
In this case $bc\le (\frac{b+c}{2})^2<(\frac{2d-a}{4})^2<d^2$.
So $\phi$ is non-PPT. Rest is simple computation.
\endproof

\begin{theorem} Let $a,b,c,d \in \mathbb{R}$ with
$a,b,c \ge 0$ and $a+b+c = 1.$ Let $\varphi =
\rho[a,b,c,d],$ $\psi= \tau [a,b,c,d].$ Let $A = \frac{1}{3}C_{\varphi}$ and
$B=\frac{1}{3}C_{\psi}$ be the Choi-Jamiolkowski matrices of $\varphi$ and $\psi$
respectively.
\begin{enumerate}[\rm(i)]
\item If $d=0$ then $A=B$ is a separable state.
\item Let $d > 0$.
\begin{enumerate}[\rm(a)]
 \item If $a \ge d$ but $bc < d^2,$ then $A$ is a non-PPT state.
 \item If $bc \ge d^2$ but $a < d,$ then $B$ is a non-PPT state.
 \item If $a \ge d,$ $bc \ge d^2$ then $A$ and $B$ are PPT states.
 \item If $a \ge d,$ $bc \ge d^2$ but $a+b < 2d$ or $a+c < 2d$ then $A$ and $B$ are PPTES.

In particular, this is true if for an arbitrary $0 < \beta <1,$ we take
$\lambda =  \beta / (\beta^2 + \beta +1),$ $a =\lambda
=d,$ $b = \lambda \beta,$ $c = \frac{\lambda}{\beta}.$

 \item If $d \le a < 2d$ and $2 (b+c) < 2d-a$ then $A$ has Schmidt number
 3 and it is non-PPT. This is equivalent to requiring:

$\frac{2}{3}<a\le 1$ together with $1-\frac{a}{2}<d\le a$, or, equivalently, $\frac{1}{2}<d\le 1$ together with $d\le a\le 1$ for
$d>\frac{2}{3}$, whereas, $2(1-d)<a\le 1$ for $\frac{1}{2}<d\le \frac{2}{3}$, and, then taking $b$ and $c\ge 0$ with $b+c=1-a$.
\end{enumerate}
\item Let $d < 0.$
\begin{enumerate}[\rm(a)]
\item If $a \ge - 2d = 2 |d|$ but $bc < d^2$ then $A$ is a non-PPT state.
\item If $bc \ge d^2$ but $a < 2 |d|,$ then $B$ is a non-PPT state.
\item If $a \ge 2 |d|,$ $bc \ge d^2,$ then $A$ and $B$ are PPT
states.
\end{enumerate}
\end{enumerate}
\end{theorem}

\proof This follows immediately from Theorem 2.6 above.\endproof

\begin{remark}\rm
This concerns the range and the rank of
$C_{{\rho}[a,b,c,d]}$ and $C_{{\tau}[a,b,c,d]}$.

(i) We first note that
\begin{align*}
&C_{\rho[a,b,c,d]} =\begin{bmatrix} a &
0&0&0&d&0&0&0&d \\[-2pt]0&c&0&0&0&0&0&0&0 \\[-2pt] 0&0&b&0&0&0&0&0&0 \\[-2pt]
0&0&0&b&0&0&0&0&0 \\[-2pt]d&0&0&0&a&0&0&0&d\\[-2pt] 0&0&0&0&0&c&0&0&0
\\[-2pt]0&0&0&0&0&0&c&0&0 \\[-2pt]0&0&0&0&0&0&0&b&0 \\[-2pt]d&0&0&0&d&0&0&0&a \\[-2pt]
\end{bmatrix},\\
&C_{\tau[a,b,c,d]} =\begin{bmatrix}
a&0&0&0&0&0&0&0&0 \\[-2pt]0&c&0&d&0&0&0&0&0 \\[-2pt]0&0&b&0&0&0&d&0&0
\\[-2pt]0&d&0&b&0&0&0&0&0 \\[-2pt] 0&0&0&0&a&0&0&0&0 \\[-2pt]0&0&0&0&0&c&0&d&0
\\[-2pt]0&0&d&0&0&0&c&0&0 \\[-2pt]0&0&0&0&0&d&0&b&0 \\[-2pt]0&0&0&0&0&0&0&0&a
\end{bmatrix}.
\end{align*}
\newpage
(ii) By (i),  the range of $C_{\rho[a,b,c,d]}$ is the linear span of
\begin{align*}
&ae_1 + de_5 + de_9,\quad ce_2,\quad be_3,\quad be_4, \\
&de_1 + ae_5 + de_9,\quad ce_6,\quad ce_7,\quad be_8 \quad \mbox{and} \\
&de_1 + de_5 + ae_9.
\end{align*}
Here $e_j$, $1\le j\le 9$ is the arrangement of product vectors
$e_p^3\otimes e_q^3$ in the lexicographic order, $(e_p^3:p=1,2,3)$ being the
standard ordered basis of $\mathbb{C}^3$.

Now the matrix $\begin{bmatrix} a & d & d \\ d & a & d \\ d
& d & a \end{bmatrix}$ has determinant $= (a-d)^2 (a+2d).$

Therefore $C_{\rho[a,b,c,d]}$ attains all ranks from $0$ to $9$
depending on the values of $a,b,c,d$ in $\mathbb{C}.$ For instance,
for $a,b,c$ all non-zero, the matrix $C_{\rho[a,b,c,d]}$ has
rank 7 if $a=d,$ $\mbox{rank}\, 8$ if $a=-2d$ and rank 9 if $a \neq
d$ and $a \neq -2d.$

(iii) Next, the range of $C_{\tau[a,b,c,d]}$ is the linear span of
$ae_1,$
$ce_2+de_4,$ $be_3+de_7,$ $de_2+be_4,$ $a e_5,$ $ce_6+de_8,$
$de_3+ce_7,$ $de_6+be_8$ and $ae_9.$ The matrices $
\begin{bmatrix} c& d \\ d & b \end{bmatrix}$ and $
\begin{bmatrix} b& d \\ d & c \end{bmatrix}$ both have
determinant $= bc-d^2.$ So the matrix $C_{\tau[a,b,c,d]}$ has rank 9
if $a \neq 0 \neq bc - d^2,$ rank 6 if $a \neq 0 = bc - d^2$ and $d
\neq 0,$ rank 6 if $a \neq 0 = bc - d^2 = d$ and $b^2 + c^2 \neq 0$,
and rank 3 if $a\neq 0 = b=c=d.$
\end{remark}

\begin{remark}\rm
 Let $a,b,c\ge 0$ and $d\in\mathbb{R}$ with $d\ne 0$. Let $A=\frac{1}{3}C_{\rho[a,b,c,d]}$ and $B=\frac{1}{3}C_{\tau[a,b,c,d]}$.
 We refer to their expanded form coming from Remark~2.8 above.
\begin{enumerate}[\rm(i)]
\item In view of Item (iii)(b) or Theorem 2.7, $A$ is a density matrix if and only if $a+b+c=1$ and
 $a\ge \max\{d,-2d\}$. We now consider only this case. Then $a>0$.
 Further, from the expanded form of $A$, it is clear that it is enough to look at its
 only non-trivial sub-block $A_1=\frac{1}{3}\begin{bmatrix}a&d&d\\d&a&d\\d&d&a\end{bmatrix}$
 acting on the span of $e_1$, $e_5$ and $e_9$. Now 3A$_1$ has eigenvalues $a+2d$, $a-d$, $a-d$. Since $d\ne 0$, we have
 $a+2d\ne a-d$. Let $\xi=\frac{1}{\sqrt{3}}(e_1+e_5+e_9)$.
 Then $A_1=\frac{1}{3}(a+2d)P_\xi+\frac{1}{3}(a-d)P$, $P_\xi$ is the projection determined by $\xi$ and $P$, the (orthogonal) projection on
 $\xi^\perp$.
 Also $\xi^\perp$ is the linear span of $e_1-e_5$ and $e_5-e_9$ and it contains no non-zero product vectors.

\item By Item (iii)(b) or Theorem 2.7, we have that $B$ is a density matrix if and only if $a+b+c=1$ and $bc\ge d^2$.
 We now consider only this case. Then $b>0$, $c>0$. Further, it is clear from
 the expanded form of $B$ that it is enough to consider its non-trivial sub-blocks $B_1=\frac{1}{3}\begin{bmatrix}c&d\\d&b\end{bmatrix}$,
 $B_2=\frac{1}{3}\begin{bmatrix}b&d\\d&c\end{bmatrix}$, $B_3=\frac{1}{3}\begin{bmatrix}c&d\\d&b\end{bmatrix}$, acting on linear spans
 $L_1,L_2,L_3$ respectively of pairs of product vectors $(e_2,e_4)$, $(e_3,e_7)$ and $(e_6,e_8)$ respectively. Because
 $d\ne 0$, $B_1$ has distinct eigenvalues. Further, any corresponding
eigenvector has the Schmidt rank~2.\ The same applies to $B_2$ and $B_3$ as well.
 \end{enumerate}
\end{remark}

\section{Power Symmetric Matrices and Construction of PPT States}
This section has been motivated by the work on ``Power symmetric matrics''
by Sinkhorn \cite{sin} and its generalization by Bapat, Jain and Prasad \cite{bab} on one hand and preservation
of positivity of a block matrix under taking powers of blocks,
the so-called Schur or Hadamard product by Choudhury \cite{ch} and Guillot, Khare and Rajaratnam \cite{gui} on
the other. The purpose is to illustrate the interesting interplay rather than the utmost generality.

\subsection{Power symmetric matrices}
We first recall the known concepts and results to be used.

\def\NN{\mathbb{N}}
\def\TT{\mathbb{T}}
\begin{enumerate}[(i)]
\item
Sinkhorn \cite{sin} calls a stochastic matrix $A$ \emph{power symmetric} if its transpose $A^t$
equals $A^q$ for some $q\in \mathbb{N}$.\ The smallest such positive integer $q$ for which this is true is called the symmetric order of $A$.

Let $A$ be a power symmetric matrix of symmetric order $q>1$.
\begin{enumerate}[(a)]
\item Sinkhorn notes that such an $A$ commutes with $A^t$ and, is therefore, normal.
\item He proves that $A$ is bistochastic by showing that a normal stochastic matrix is bistochastic.
\item He also shows that $A$ satisfies $A^{q^2}=A$, and, therefore, an eigenvalue $\lambda$ of $A$ is either zero or has modulus one.
\item For $r\in \NN$, let $J_r$ be the bistochastic $r\times r$ matrix having all entries equal to $\frac{1}{r}$.
Sinkhorn obtains a general form for $A$ in terms of certain $J_r$'s.
\end{enumerate}

\item It may be of interest to note that Sinkhorn is famous for his Diagonal Theorems:
If $A$ is an $n\times n$ matrix with strictly positive elements then there
exist diagonal matrices $D_1$ and $D_2$ with strictly positive elements
such that $D_1A D_2$ is doubly stochastic, $D_1$ and $D_2$ are essentially
unique i.e., up to a positive multiple by $a$ and $\frac{1}{a}$ respectively (Ann. Math. Statistics, 1964).\
Further, there exists a unique stochastic matrix
of the form $DAD$, where $D$ is a diagonal matrix with positive diagonal elements (Can. J. Math., 1966).

This stimulated further research by reputed mathematicians like R. Brualdi.

\item The property (i)(a) persists even if we take $A$ to be a real matrix instead of stochastic. Further, (i)(c) holds even if
we take $A$ to be a complex matrix satisfying $A^t=A^q$ or $A^*=A^q$ for some $q\in \NN$ with $q>1$.

This motivates our next result.

For $q\in \NN$ with $q>1$, let $\TT_q^{(1)}=\{\lambda:\lambda^{q-1}=1\}$ and
$\TT_q^{(2)}=\{\lambda: \lambda^{q^2-1}=1\neq \lambda^{q-1}\}$. We note that
$\TT_q^{(2)}$ can be expressed as a disjoint union of
pairs $T_r=\{\lambda_r,\lambda_r^q\}$, $1\le r\le d_q$, say. A subset $S_2$ of
$\TT_q^{(2)}$ will be called \emph{unpaired} if $\# S_2\cap T_r\le 1$
for $1\le r\le d_q$.

For the sake of convenience an empty sum will be taken to be zero. Also for
an $A=[a_{jk}]$ in $M_n$, $\bar A$ denotes the matrix $[\bar a_{jk}]$.
\end{enumerate}

\begin{theorem}
Let $A$ ba a non-zero real or normal
matrix in $M_n$ and $q\in \NN$ with $n,q>1$. Then $A$ satisfies $A^t=A^q$ if and only if there exist
subsets $S_1$ and $S_2$ of $\TT_q^{(1)}$ and $\TT_q^{(2)}$ respectively and mutually orthogonal non-zero projections $P_s$,
$s\in S_1\cup S_2$ such that
\begin{enumerate}[\rm(i)]
\item $S_2$ is unpaired;
\item for $s\in S_1$, $P_s$ is a real matrix;
\item for $s\in S_2$, $\bar P_s$ and $P_{s'}$ are orthogonal for $s'\in S_1\cup S_2$;
\item $A=\sum\limits_{s\in S_1} sP_s+\sum\limits_{s\in S_2} sP_s+\sum\limits_{s\in S_2}s^q \bar P_s$.
\end{enumerate}
\end{theorem}

\begin{proof}
The `if' part follows immediately by taking the transpose of $A$ as in (iv).

For the other part, suppose $A^t=A^q$. Then by 3.A(iii) $A$ is normal and $A^{q^2}=A$. Let $\sum =$ the set of non-zero eigenvalues of
$A$. Because $A\neq 0$, $\sum$ is nonempty. Also $A^{q^2}=A$ forces $\sum$ to be contained
in $\{\lambda:\lambda^{q^2-1}=1\}=\TT_q^{(1)}\cup \TT_q^{(2)}$.
Let $S_1=\sum \cap \TT_q^{(1)}$ and $\sum_2 =\sum \cap \TT_q^{(2)}$. Because $A$ is normal,
there exist non-zero mutually orthogonal projections
$\{P_s:s\in \sum\}$ such that $A=\sum\limits_{s\in \sum}sP_s$. Then $A^t =\sum\limits_{s\in \sum}sP_s^t=\sum\limits_{s\in \sum}s \bar P_s$.
Also $A^q=\sum\limits_{s\in \sum}s^q P_s$.

But $A^t=A^q$. So $\sum\limits_{s\in \sum}s \bar P_s=\sum\limits_{s\in \sum}s^q P_s$.

So for each $s\in \sum$, there exists $\lambda_s\in \sum$ satisfying $s=\lambda^q_s$ and $\bar P_s=P_{\lambda_s}$,

For $s\in S_1$, $\lambda_s=s$ and, therefore, $\bar P_s=P_s$, i.e., $P_s$ is a real matrix.

For $s\in \sum_2$, $\bar P_s=P_{\lambda_s}$; and also, $\lambda_s\in \sum_2$ with $\bar P_{\lambda_s}=P_s$.
We can form an unpaired set $S_2\subset \sum_2$ with $\{s,s^q: s\in S_2\}=\sum_2$.

Hence $A=\sum\limits_{s\in S_1}sP_s+\sum\limits_{s\in S_2}sP_s+\sum\limits_{s\in S_2}s^q\bar P_s$ and conditions (i), (ii), (iii) are also satisfied.
\end{proof}

\begin{definition}\rm
Let $A$ be a normal matrix. For a $q\in \NN$, $q>1$, $A$ will be called \emph{$q$-reflected} if $A^t=A^q$ and it will be
called \emph{reflected} if it is $q$-reflected for some $q\ge 2$.
\end{definition}

\begin{remark}\rm
Let $A$ be a non-zero matrix..
\begin{enumerate}[\rm(i)]
\item If $A$ is power symmetric of symmetric order $>1$ then it is reflected.
\item If $A$ is reflected then $A$ is a partial isometry.
\item Suppose $A$ is reflected. Then  the following hold.
\begin{enumerate}[(a)]
\item $A$ is positive if and only if $A$ is projection and real.
\item $A$ is Hermitian if and only if $A=P-Q$ for some real mutually orthogonal projections $P$ and $Q$.
\item If $A$ is $q$-reflected, then, in the notation of Theorem 3.1,  A satisfies (i) to (iv) of Theorem~3.1 and  also,
\begin{align*}
\bar AA
&=\sum_{s\in S_1}P_s+\sum_{s\in S_2} s\bar s^qP_s+\sum_{s\in S_2}\bar s s^q \bar P_s=A\bar A\\
&=\sum_{s\in S_1}P_s+\sum_{s\in S_2}\bar s^{q-1}P_s+\sum_{s\in S_2} s^{q-1}\bar P_s~;
\end{align*}
in particular, $\bar A A$ is a real $q$-reflected matrix.
\end{enumerate}
\end{enumerate}
\end{remark}

\subsection{Generalized power symmetric matrices}

We begin with the relevant material from Bapat, Jain and Prasad  \cite{bab}.

\begin{enumerate}[(i)]
\item As defined by Bapat, Jain and Prasad \cite{bab}, a
\emph{generalized power symmetric} matrix is a stochastic matrix $A$ which satisfies
$(A^p)^t =A^q$ for some $p,q\in \NN$ with $p<q$.
\item Bapat et al \cite{bab} obtain several properties of a generalized
power symmetric matrix $A$ and also give a general form for $A$
in terms of $J_r$'s.\
An easy one to be used is, $A^{p^2}=A^{q^2}$ and therefore, an eigenvalue $\lambda$ of $A$
has to be either zero or satisfy $\lambda^{q^2-p^2}=1$.
We may formulate an analogue of Theorem~3.1 but it is not so neat or revealing.
\end{enumerate}

We make an attempt to utilize contents of this section to construct PPT states.\
For this purpose we give methods to produce classes
of positive block matrices that are positive under partial
transpose, in short, PPT.\
The idea is to replace transpose by a power in a broad sense and see
if that will be helpful.

\subsection{Positivity of block matrices and Schur or Hadamard products}

We collect a few results in this direction.
\begin{enumerate}[(i)]
\item This is a well-known result and one may find it in books like \cite{rbh} and \cite{hia}:

Let $A,B,C$ be $n\times n$ matrices. Then the block matrix $R=\begin{bmatrix}A&B\\B^*&C\end{bmatrix}\ge 0$
if and only if $A\ge 0$, $C\ge 0$ and there exists a contraction $K$ such that $B=A^\frac{1}{2}KC^\frac{1}{2}$.

\item \textsc{Choudhury's Theorem.}\ Choudhury \cite{ch} considered a block matrix $H=[H_{jk}]$, where $H_{jk}$'s are normal $n\times n$ matrices
for $1\le j$, $k\le m$.\ Theorem~5 \cite{ch} says that if the $m^2$
matrices $\{H_{jk}:1\le j,k\le m\}$ are a commuting family and $H$ is positive
semi-definite then
$H_\alpha=[H^\alpha_{jk}]$ is positive semi-definite for all $\alpha=1,2,\ldots$.

\item
Arbitrary Hadamard powers can be defined for non-negative matrices (cf. [4]). In that spirit, Guillot,
Khare and Rajaratanam \cite{gui} give interesting further developments, but reveal that
actions of replacing
$\alpha$ by other numbers or relaxing the conditions limit the possibility of preserving
the positive semi-definiteness and that of its application to the
construction of PPT states.
\end{enumerate}

\def\cM{\mathcal{M}}
\def\cH{\mathcal{H}}
\begin{definition}\rm
Let $\cM$ be a set of normal $n\times n$ matrices.
\begin{enumerate}[\rm(i)]
\item For $q\in\NN$, $q>1$, the set $\cM$ will be called \emph{$q$-reflected} if $A^t=A^q$ for each $A$ in $\cM$.
It will be called \emph{reflected} if it is {$q$-reflected} for some $q\ge 2$.

\item For $p,q\in\NN$ with $p<q$, the set $\cM$
will be called \emph{generalized $(p,q)$-reflected} if $(A^p)^t =A^q$ for each $A$ in $\cM$.\
Further it will be called \emph{generalized reflected} if it
is $(p,q)$-reflected for some $p,q\in \NN$ with $1\le p<q$.
\end{enumerate}
\end{definition}

\begin{theorem}
Let $\cH=\{H_{jk}:1\le j, k\le m\}$ be a reflected commuting family of matrices such that $H=[H_{jk}]$ is positive. Then $H$ is positive under partial transpose.
\end{theorem}

\proof
There exists a $q\in\NN$ with $q\ge 2$
such that $\cH$ is $q$-reflected.\ We first note that
it follows from Remark~3.3(iii)(a)
that each $H_{jj}$ is real and also a projection.\
So $H_{jj}^t=H_{jj}=H_{jj}^q$ for $1\le j\le m$. By Choudhury's Theorem~3.C(ii) $H_q$ is positive.
But $H_q=H^{PT}$. So $H$ is positive under partial transpose.
\endproof

\begin{theorem}
Let $B$ be reflected and $A,C\ge 0$ satisfy $R=\begin{bmatrix}A&B\\ B^*&C\end{bmatrix}\ge 0$ and $B$ commutes with $A$ or $C$.
If either $B$ or both $A$ and $C$ are real, then $R$ is positive under partial transpose.
\end{theorem}

\begin{proof}
By 3.C(i) there exists a contraction $K$ such that $B=A^\frac{1}{2}KC^\frac{1}{2}$. Suppose $B$
commutes with $A$.\ Then for $q\in\NN$, $B^q=A^\frac{1}{2}B^{q-1}KC^\frac{1}{2}$.
By Theorem~3.1, $B$ is a contraction.\ So $B^{q-1}K$ is a contraction.\
By 3.C(i), $V=\begin{bmatrix}A&B^q\\ (B^q)^*&C\end{bmatrix}\ge 0$.

Since $B$ is reflected, $B^t=B^q$ for some $q\in \NN$.
So $\begin{bmatrix}A&B^t\\ (B^t)^*&C\end{bmatrix}=V\ge 0$.
Now $R^{PT}$ is $V$ in case both $A$ and $C$ are real, and $\bar V$ if $B$ is real.
Hence $R$ is positive under partial transpose.

Similar arguments give the result when $B$ commutes with $C$.
\end{proof}

\begin{remark}\rm
Given any $B\in M_n$, there exist real $A,C\ge 0$ with $B$ commuting with $A$ and $C$ such that
$\begin{bmatrix}A&B\\ B^*&C\end{bmatrix}\ge 0$. All we have to do is to note that
$R_B=\begin{bmatrix}0&B\\B^*&0\end{bmatrix}$ is Hermitian and therefore there exists $\lambda\ge 0$ such that $-\lambda I_{2n}\le R_B\le \lambda I_{2n}$.
We can take $A=C=\lambda I_n$.
\end{remark}

\begin{theorem}
Let $A,B,C\in M_n$ be such that $B$ is generalized reflected and invertible, $A,C\ge 0$ satisfy
$R= \begin{bmatrix}A&B\\ B^*&C\end{bmatrix}\ge 0$. If $B$ and $B^t$ commute with $A$ or $C$
and also either $B$ or both $A$ and $C$ are real, then $R$ is positive under partial transpose.
\end{theorem}

\begin{proof}
Because $B$ is generalized reflected there exist $p,q\in\NN$ with $p<q$
such that $(B^p)^t=B^q$.
Properties for generalized power symmetric matrices noted in 3.B(ii) hold for $B$.
To elaborate, eigenvalues $\lambda$ of $B$ satisfy $\lambda^{p^2}=\lambda^{q^2}$.\
Because $B$ is normal this gives $\|B\|=1$.
Now $B^{-1}$ is generalized reflected as well. So $\|B^{-1}\|=1$ too.

By 3.C(i), the condition $R\ge 0$ gives a contraction $K$ such that $B=A^\frac{1}{2}KC^\frac{1}{2}$.
Now $(B^p)^t=B^q$, So $B^t=((B^{p-1})^t)^{-1} B^q=((B^t)^{-1})^{p-1}B^q$.
Suppose $B$, $B^t$ commute with $A$.
Then
\begin{align*}
B^t
&=((B^t)^{-1})^{p-1} B^{q-1}B=((B^t)^{-1})^{p-1} B^{q-1} A^\frac{1}{2}KC^\frac{1}{2}\\
&=((B^t)^{-1})^{p-1} A^\frac{1}{2} B^{q-1} KC^\frac{1}{2}\\
&=A^\frac{1}{2} ((B^t)^{-1})^{p-1} B^{q-1} KC^\frac{1}{2}\,.
\end{align*}
As already noted $B$ and $B^{-1}$ are contractions. So,
$K_1=((B^t)^{-1})^{p-1} B^{q-1}K$ is a contraction. By 3.C(i),
$V=\begin{bmatrix}A&B^t\\ (B^t)^*&C\end{bmatrix}\ge 0$.

Hence as argued in the proof Theorem 3.6 above $R$ is positive under partial transpose.
\end{proof}

\begin{remark}\rm
We can relax the condition of invertibility on $B$ and work with
the group inverse $G$ of $B$.\
\end{remark}

\begin{remark}\rm
Let $A,B,C,R$ be as in 3.C(i) with $\tr R$ not equal to zero.
H.J. Woerdeman \cite{w} gives conditions under which $\rho=\frac{1}{\tr R}\,R$ is a separable state.
We are now in a position to use them, particularly \cite{w}, Theorem~3.2 onwards.  We record the simplest case:

Let $A=I, B, C, R$ be as in Theorem~3.6 (Theorem~3.8). Then $\rho$ is separable.\ To see this, let
$W=\begin{bmatrix}
I& B^*\\B &C\end{bmatrix}$ and
$V=\begin{bmatrix}
I& B^t\\(B^t)^*&C\end{bmatrix}$.\
Then $W=V$ if $B$ is real and $\bar V$ if $C$ is real.\
So as in the proof of Theorem~3.6 or 3.8 we have $W\ge 0$. We can now
apply the discussion in \cite{w}.
\end{remark}

\section{Peres Condition and Unitary Equivalence of a Matrix to its Transpose}

We follow the notation and terminology of Garcia and
Tener \cite{gt} who obtained a canonical decomposition for complex
matrices $T$ which are \emph{UET}, i.e., \emph{unitarily equivalent to their transpose} $T^t
(UET).$

\subsection{}
We collect a few facts from \cite{gt} for ready reference.
\begin{enumerate}[\rm(i)]
\item \cite[\S1]{gt}. In his problem book (\cite{halmos}, Pr. 159)
Halmos asks whether every square matrix
is UET and in his discussion gives the counterexample $
\begin{pmatrix} 0 & 1 & 0 \\ 0 & 0 & 2 \\ 0 & 0 & 0 \end{pmatrix}$,
which is not UET. Every Toeplitz matrix is UET via the
permutation matrix which reverses the order of the standard basic
vectors.

\item {\cite[Theorem 1.1]{gt}}. A matrix $T$ in $M_n$ is UET if and only
if it is unitarily equivalent to a direct sum of (some of the summands may be
absent):
\begin{enumerate}[\rm(a)]
\item irreducible complex symmetric matrices (CSMs),
\item irreducible skew-Hamiltonian matrices (SHMs) (such matrices are necessarily
$8 \times 8$ or larger, a SHM is a $2d \times 2d$ block matrix of
the form $\begin{pmatrix} A & B \\ D & A^{t} \end{pmatrix}$ with $B^t = - B$ and $D^t = -D),$

\item $2d \times 2d$ blocks of the form $\begin{pmatrix} A & 0
\\ 0 & A^{t} \end{pmatrix}$ where $A$ is irreducible and
neither unitarily equivalent to a complex symmetric matrix (UECSM)
nor unitarily equivalent to a skew-Hamiltonian matrix
(UESHM) (such matrices are necessarily $6 \times 6$ or larger).
\end{enumerate}

Moreover, the
unitary orbits of the three classes described above are pairwise
disjoint.

\item (\cite{gt}, Corollary 2.3).\ If $T$ is UET and has order $n
\times
n$ with $n \leq 7,$ then $T$ is UECSM.

\item (\cite{gt}, 8.3 and 8.4).\ $S$ is UET if and only if $S$ is
unitarily equivalent to a matrix $T$ that satisfies $TQ=QT^t,$ where
$Q$ is a unitary matrix of the special form (some of the blocks may
be absent and empty blocks are all zero):
$$
Q = \left (\begin{array}{c|c|c|c|c} Q_{+} & & & & \\ \hline &Q_{-} & & & \\
\hline
& & \begin{array}{cc}0 & \lambda_1 X_1^t \\ X_1 & 0
\end{array}& & \\ \hline
& & & \ddots& \\ \hline & & & &
 \begin{array}{ccc} 0 & \lambda_r X_r^t \\ X_r & 0
\end{array} \end{array} \right ),$$
where
\begin{enumerate}[\rm(a)]
\item $Q_{+} = Q_{+}^t$ is complex symmetric and unitary,
\item $Q_{-} = - Q_{-}^t$ is skewsymmetric and unitary,
\item $\lambda_i \neq \pm 1$ and $X_i$ is unitary for $i =1,2, \ldots, r. $
\end{enumerate}
\newpage
\item (\cite{gt}, 8.5). Given $Q$ as in (iv) above, $T$ is as in (iv)
above if and only if
 $$T = \left (\begin{array}{c|c|c|c|c} T_{+} & & & & \\ &T_{-} & & & \\ \hline
& & \begin{array}{cc} A_1 & 0 \\ 0 & X_1 A_1^t X_1^{\ast}
\end{array} & & \\ \hline & & & \ddots& \\ \hline & & & &
\begin{array}{ccc} A_r & 0 \\ 0 & X_r A_r^t X_r^{\ast}
\end{array} \end{array} \right ),$$
where
\begin{enumerate}[\rm(a)]
\item $T_{+} = Q_{+} T_{+}^{t} Q_{+}^{\ast}$ (such a $T_{+}$ is UECSM),
\item $T_{-} = Q_{-} T_{-}^t Q_{-}^{\ast}$ (such a $T_{-}$ is UESHM),
\item $A_1, \ldots, A_r$ are arbitrary.
\end{enumerate}

In fact this is the final step of the proof of (ii) in \cite{gt}.
\end{enumerate}

It is my pleasure to thank my students Priyanka Grover and Tanvi
Jain. Priyanka came to discuss UET in some other context and Tanvi
found a former version of \cite{gt} on the internet for that context.

\begin{definition} \rm
A tuple $(Y_1, \ldots, Y_s)$ of $n
\times n$ matrices is said to be {\it collectively unitarily
equivalent to the respective transposes (CUET)} if there is a unitary
$U$ with $Y_j = U Y_j^t U^{\ast}$ for $ 1 \leq j \leq s.$
\end{definition}

\begin{remark}\rm
W.B. Arveson \cite[Lemma
A.3.4]{wba} gives that $\left(\!\!\begin{pmatrix}0 & \lambda &
1 \\ 0 & 0 & 0 \\ 0 & 0 & 0 \end{pmatrix}\!,\!
\begin{pmatrix} 0&0& \mu \\ 0 & 1 & 0 \\ 0 & -\lambda & 0
\end{pmatrix}\!\!\right )$ is not CUET where $\lambda$ is a
non-real
complex number and $\mu$ is a complex number with $|\mu | = (1 +
|\lambda|^2)^{\frac{1}{2}}.$
\end{remark}

\begin{remark}\rm
Discussion in 4.A above tells us how to construct CUET tuples viz.,
choose a $Q$ as in item 4.A(iv) and then $Y_j$'s as in
4.A(v) by varying $T_{+}, T_{-}, A_k$ for $1 \leq
k \leq r.$
\end{remark}

\begin{theorem} Let $[A_{jk}]$ be a positive block
matrix such that $(A_{jk}: 1 \leq j,k \leq n)$ is CUET. Then
$[A_{jk}^t]$ is positive.
\end{theorem}

\proof There is a unitary matrix $U$ such that
$A_{jk} = UA_{jk}^t U^{\ast}$ for $1 \leq j,k \leq n.$ Let
$\widetilde{U}$ be the block matrix $[\delta_{jk}
U],$ with $\delta_{jk} = 0$ for $j \neq k$ and $1$ for
$j=k.$ Then $\widetilde{U}$ is unitary and $[A_{jk}^t] =
\widetilde{U}^{\ast} [A_{jk}] \widetilde{U}.$ So
$[A_{jk}^t]$ is positive.
\endproof

\begin{construction}\rm Remark 4.3 and Theorem
4.4 put together tell us how to construct PPT matrices.
\begin{enumerate}[\rm\emph{Step}~1:]
\item Let $n \ge 2$ and put $m=\frac{n(n+1)}{2}.$
Use Remark 4.3 to construct a CUET $m$-tuple $(Y_j:1 \leq j
\leq m)$ with $Y_j \ge 0$ for $1 \leq j \leq n$ (we may take all
$Y_j,$ $1 \leq j \leq n$, to be zero, for instance). For $1 \leq j
\leq m,$ we have $Y_j = QY_j^t Q^{\ast}$ and, therefore $Y_j^{\ast} =
QY_j^{\ast t} Q^{\ast}.$ We set $B_{jj} = Y_j$ for $1 \leq j \leq n,$
arrange $Y_j$ for $n+1 \leq j \leq \frac{n(n+1)}{2}$ as $B_{pq},$ $1
\leq p < q \leq n$ and take $B_{qp} = B_{pq}^{\ast}$ for $1 \leq p <
q \leq n.$ Thus, we obtain a block matrix $B = [B_{jk}]$ which is
Hermitian and $\{B_{jk} : 1 \leq j, k \leq n \}$ is CUET.

\item The set $\{a
\in \mathbb{R}: B+a \,\, I_{n^{2}} \ge 0\}$ is an interval $[a_0,
\infty)$ for some $a_0 \in \mathbb{R}.$ We take any $a$ in this
interval and set $A=B+a \,\, I_{n^{2}}$ i.e., $A_{jk} = B_{jk}$ for $j
\neq k,$ whereas $A_{jj} = B_{jj} + a I_n$ for $1 \leq j, k \leq n.$
Then $\{A_{jk} : 1 \leq j, k \leq n \}$ is CUET and $A \ge 0.$ So we
can apply Theorem 4.4 to conclude that $A$ is a PPT matrix.
\end{enumerate}
\end{construction}

\begin{remark}\rm
We can, of course, formulate ``unitarily equivalent'' versions of various notions in \S3 above for
subsets $\cM$ of $M_n$. 
Suitable variants of results in \S3 can then be obtained to give new classes of positive matrices with positive
partial transpose and thus, PPT states. We do not go into the details here.
\end{remark}

\section{Quantum Dynamical Semigroups involving Separable and Entangled States}

Let $\mathcal{H}$ be a Hilbert space and $\tau$ the transpose map on
$\mathcal{B}(\mathcal{H})$ with respect to some orthonormal basis for
$\mathcal{H}.$ Let $\ast$ or $\dagger$ be the adjoint map on
$\mathcal{B}(\mathcal{H})$ that takes $x$ to $x^{\ast}.$ Let $\mathcal{X}$ be a
linear subspace of $\mathcal{B}(\mathcal{H})$ which is closed under $\tau$ as
well as $\ast.$ We shall consider $C_0$-semigroups $(T_t)_{t \ge 0}$
as well as $T_0$-constricted $C_0$-semigroups $(T_t)_{t \ge 0}$
of operators on $\mathcal{X}$ to itself.

We begin with a few examples.

\begin{examples}\rm
\begin{enumerate}[\rm(i)]
\item This is modelled on St{\o}rmer's Example 8.13 \cite{es}
and is in a foliated
form with $\Lambda_{t1} = \begin{bmatrix} 1 & 0 \\ 1 -
e^{-t} &
e^{-t} \end{bmatrix}$, $\Lambda_{t2} = e^{-\frac{t}{2}}
\mathcal{I}_{F_{2}}$ and $\Lambda_t = \Lambda_{t1} \oplus
\Lambda_{t2}$ in the notation of item 2.H(i)(b). It is a non-PPT quantum
dynamical semigroup.

\item If we are interested in separable maps we have to do away
with the condition
$T_0 = \mathcal{I}d,$ which we now do.
\end{enumerate}

This example is modelled on the example of the two spin
$\frac{1}{2}$-states given by Horodecki et al
\cite{hhh}. It is in a foliated form with $\Lambda_1^{p,a,b} =
\begin{bmatrix} pa^2 & (1-p)b^2 \\ (1-p) a^2 & pb^2 \end{bmatrix}$ and
$\Lambda_2^{p,a,b} = \begin{bmatrix}
pab & (1-p)ab \\ (1-p) ab & pab \end{bmatrix}$ with $0
\leq p \leq 1,$ $a > 0,$ $b > 0$ and $\Lambda^{p,a,b} = \Lambda_1^{p,a,b}
\oplus\Lambda_2^{p,a,b}.$
\end{examples}

\newpage
Taking Pauli matrices $\sigma_0 = I_2,$ $\sigma_1 = \begin{pmatrix} 0
& 1 \\ 1 & 0 \end{pmatrix},$ $\sigma_2 = \begin{pmatrix} 0
& -i \\ i & 0 \end{pmatrix},$ $\sigma_3 = \begin{pmatrix} 1
& 0 \\ 0 & -1 \end{pmatrix}$ as a basis, the map
$\Lambda^{p,a,b}$ has
the simple form
$$\begin{pmatrix} \frac{a^2 + b^2}{2} & 0 & 0 & \frac{a^2- b^2}{2}
\\ 0 & ab & 0 & 0 \\ 0 & 0 & 2 (\frac{1}{2}-p)ab & 0 \\ (p - \frac{1}{2})(a^2 -
b^2) & 0 & 0 &(p - \frac{1}{2})(a^2 + b^2) \end{pmatrix}. $$

(a) The map $\wedge^{p,a,b}$ is unital if and only if $a^2+b^2=2$ and
either $p=\frac{1}{2}$ or $a^2=b^2$ $(=1)$. On the other hand, the map is
trace preserving if and only if $a^2=1=b^2$.

(b) As noted by Horodecki et al, it is a separable map if and only if $p =
\frac{1}{2}$
and, in that case, the matrix becomes $\begin{pmatrix}
\frac{a^2 + b^2
}{2} & 0 & 0 & \frac{a^2 - b^2}{2} \\ 0 & ab & 0 & 0 \\ 0&0&0&0 \\ 0&0&0&0
\end{pmatrix}$ and thus the range is the commutative algebra spanned by
$\sigma_0$ and $\sigma_1.$

(c) Taking $a^2+b^2=1,$ $a = \cos \theta,$ $b = \sin
\theta,$ $0 < \theta \leq \frac{\pi}{4},$ $u = \sin 2 \theta$ motivates
the semigroup
$$T_t = \left ( \frac{1}{2} \right )^{t}
\begin{bmatrix}1 & 0 & 0 &\sqrt{1-u^2} \\ 0 & u^t & 0 & 0 \\
0&0&0&0 \\ 0&0&0&0 \end{bmatrix}, \quad t \ge 0.$$
We note that $T_0$ is the idempotent
$$
\begin{bmatrix} 1 & 0 & 0 &\sqrt{1-u^2} \\ 0 & 1 & 0 & 0 \\
0&0&0&0 \\ 0&0&0&0 \end{bmatrix}.$$
(d) We may consider the variant for $t \ge 0,$
$$
S_t = \begin{bmatrix} 1 & 0 & 0 &\sqrt{1-u^2} \\ 0 & u^t & 0 & 0 \\
0&0&0&0 \\ 0&0&0&0 \end{bmatrix}.
$$
All these take $I_2$ to itself and also $S_0 = T_0.$
Any one of $S_t$'s is (and thus, all are) trace preserving if and only if $u=1$,

i.e., $S_t$'s are all the same as the projection given by the matrix
\begin{align*}
\begin{bmatrix}
1&0&0&0\\
0&1&0&0\\
0&0&0&0\\
0&0&0&0
\end{bmatrix}.
\end{align*}

\begin{theorem}[Trichotomy]
Let $(T_t)_{t \ge 0}$ be a $T_0$ constricted
quantum dynamical semigroup. Then one and only one of the following holds.
\begin{enumerate}[\rm(i)]
\item For each $t \ge 0, $ $T_t$ is separable.
\item There exists $t_0 > 0$ such that $T_t$ is entangled for $t < t_0$ but $T_t$ is separable for $t \ge t_0.$
\item For each $t \ge 0,$ $T_t$ is entangled.
\end{enumerate}
Moreover, (i) holds if and only if $T_0$ is separable.
\end{theorem}

\proof The three conditions are mutually
exclusive. Suppose (iii) does not hold. Then the set $S = \{ t \ge 0:T_t$ is separable$\} \neq \phi$.
If $t \in S$ then for $s> t,$ $T_s = T_{s-t} T_{t}.$ By Theorem 2.1, $T_s$ is separable. So $S$
is an interval of the form
$(t_0, \infty)$ or $[t_0, \infty).$ By item 2.D(vi) and the
condition of strong continuity on $(T_t)_{t \ge 0},$ $T_{t_{0}}$ is
separable. So $S = [ t_0, \infty).$ If (i) does not hold, we have
$t_0 > 0.$ Thus (ii) holds.

Now suppose $T_0$ is separable. Then the set $S = \{ t \ge 0:
T_t$ is separable$\}$ contains $0$.
As seen above, if $S \neq \phi,$ then $S= [t_0, \infty)$ for some
$t_0 \in [0, \infty)$.\ So $S = [0, \infty)$ i.e. (i) holds.
\endproof

\begin{definition}\rm
The space $\mathcal{X}$ will be
said to be {\it normal} if for $x \in \mathcal{X},$ $\# \{x, x^{\ast},
xx^{\ast}, x^{\ast} x \} \cap \mathcal{X} \leq 3.$ In other words,
each $x \in \mathcal{X}$ is either normal or else at most one of
$xx^{\ast}$ and $x^{\ast}x$ is in $\mathcal{X}.$
\end{definition}

\begin{proposition}
Let $\varphi$ be an idempotent
$\ast$-map on $\mathcal{X}.$ If $\varphi$ is co-Schwarz then the
range of $\varphi$ is normal.
\end{proposition}

\proof Let $\mathcal{Y} = \varphi (\mathcal{X}).$ Since
$\varphi$ is a $\ast$-map, for $x \in \mathcal{Y},$ $x^{\ast}$ is in
$\mathcal{Y}.$ Since
$\varphi^2 = \varphi$ we have $\varphi | \mathcal{Y} = \mathcal{I}
d_{\mathcal{Y}}.$ Let, if
possible, there exist $y \in \mathcal{Y}$ with $y^{\ast} y,$
$yy^{\ast} \in \mathcal{Y}.$ Since $\varphi$ is co-Schwarz, we have
$\tau \varphi
(y^{\ast} y) \geq \tau \varphi (y^{\ast}) \tau \varphi (y)$,
i.e. $\tau (y^{\ast} y) \geq \tau(y^{\ast}) \tau (y).$
So $\tau (y) \tau (y^{\ast}) \ge \tau (y^{\ast}) \tau (y).$
We may interchange the role of $y$ and $y^{\ast}$ and get
$\tau(y^{\ast}) \tau(y) \ge \tau (y) \tau (y^{\ast}).$
So $\tau(yy^{\ast}) = \tau (y^{\ast}y).$ Therefore $yy^{\ast} =
y^{\ast} y.$
\endproof

\begin{remark}\rm Let $(T_t)_{t \ge 0}$ be a
$T_0$-constricted quantum dynamical semigroup.
\begin{enumerate}[\rm(i)]
\item If the range $\mathcal{R}_0$ of $T_0$ is normal, then
for $t>0$, the range $\mathcal{R}_t$ of $T_t$ is normal simply because
$\mathcal{R}_t \subset \mathcal{R}_0$.

\item One can have more Trichotomy results by replacing
``separable'' by
\begin{enumerate}[\rm(a)]
\item PPT, or
\item has Schmidt number $\leq r,$ or
\item has normal range
\end{enumerate}
and then ``entangled'' by the corresponding negations like non-PPT,
has Schmidt number $> r$ and has non-normal range.

\item In fact, the first condition in any such Trichotomy
holds if and only if it holds for $T_0.$ By 2.D(viii), it holds if
$\mathcal{R}_0$ is contained in an abelian $C^{\ast}$ algebra acting
on a separable Hilbert space $\mathcal{H}.$

\item A non-commutative $C^{\ast}$-algebra is not normal. So
for an interesting theory, we can give up the condition (i) of
Trichotomy and instead take $T_0 = \mathcal{I}d.$
\end{enumerate}
\end{remark}

\begin{theorem}
Let $\mathcal{X}$ be a
non-commutative $C^{\ast}$-algebra and $(T(t))_{t \ge 0}$ be a
quantum dynamical semigroup of  completely positive maps. If
for some $t_0 > 0,$
$T(t_0)^{-1}$ exists and is a Schwarz map, then for each $t > 0,$
$T(t)$ is non-PPT.
\end{theorem}

\proof We refer to item 2.G(vii)(a) as for the
proof of Theorem~2.4. We use the fact that the product of two Schwarz maps is
a Schwarz map. For $0
< t < t_0,$ $T(t)^{-1} = T (t_0-t)(T(t_0))^{-1},$ and therefore,
$T(t)^{-1}$ is a Schwarz map.\ Also for $n \in \mathbb{N},$ $0 < s <
t_0,$ $t = nt_0 + s,$
$T(t)^{-1} = T(s)^{-1} (T(t_0)^{-1})^n,$
and therefore, $T(t)^{-1}$ is a Schwarz map. Let, if possible, for
some $t > 0,$ $T(t)$ be PPT. Then $\tau T(t)$ is completely positive.
So $\tau = \tau T(t) (T(t))^{-1}$ is a Schwarz map, which is not so
because $\mathcal{X}$ is non-commutative.
\endproof

We now illustrate results in this section with examples of
generalized Choi maps discussed in 2.H above.

\begin{example}\rm
\begin{enumerate}[\rm(i)]
\item This may be thought of as continuation of 2.H.
We begin by recalling relevant details, which are well-known from the
theory of circulant matrices (cf. \cite{rbbtesr}, \cite{rbh},
\cite{gkps}).

\item Let $\alpha \in \mathbb{R},$ $\beta \in \mathbb{C}.$ Then
\begin{align*}
D (\alpha, \beta, \bar{\beta})
&= \alpha I_3 + \beta (E_{12} +E_{23} + E_{31}) + \bar{\beta} (E_{21} +E_{32} + E_{13})\\
&= \alpha I_3 + \beta L + \bar{\beta} L^{\ast},
\end{align*}
where $L = E_{12} + E_{23} + E_{31}$.

We note that $L^2=L^*$, $LL^{\ast} =L^{\ast}L = I_3.$ So $L$ is a
unitary matrix with eigenvalues $1, \omega,\omega^2$ and is expressible as $L=W
\Diag (1,\omega,\omega^2)W^{\ast}$ with $W = \frac{1}{\sqrt{3}}
\begin{bmatrix}
 1 & 1 & 1 \\ 1 & \omega & \omega^2 \\ 1 & \omega^2 & \omega
\end{bmatrix}.$

Here $\omega=-\frac{1}{2}+\frac{\sqrt{3}}{2}i$, a cube root of unity.\\
So $D(\alpha, \beta, \bar{\beta}) = W \Diag (\alpha + \beta +
\bar{\beta}, \alpha+ \beta \omega + \bar{{\beta}}\omega^2, \alpha+\beta \omega^2 +
\bar{\beta} \omega ) W^{\ast}.$

\item For $a,b,c \in \mathbb{C},$ $D (a, b, c)^{\ast} = D
(\bar{a}, \bar{c},\bar{b})$ and thus, $D=D(a,b,c)$ is
normal.
Further, $D_1 = \frac{1}{2} (D+D^{\ast}) = D (\Re \,a, \beta,
\bar{\beta})$ with $\beta = \frac{1}{2}(b+\bar{c}),$ and $D_2 =
\frac{1}{2i} (D-D^{\ast}) = D (\mathrm{Im} \, a, \gamma, \bar{\gamma})$ with
$\gamma = \frac{1}{2i} (b-\bar{c}).$ So by (ii),
\begin{align*}
D(a,b,c) &= D_1 + iD_2\\
&= W \Diag (a+(\beta + i \gamma) + (\bar{\beta} + i \bar{\gamma}), a+ (\beta + i \gamma) \omega + (\bar{\beta} + i \bar{\gamma})\omega^2,\\
&\quad a+(\beta + i \gamma)\omega^2+(\bar{\beta} + i \bar{\gamma}) \omega)W^{\ast}\\
&= W \Diag (a+b+c, a+b\omega+c\omega^2, a+ b\omega^2+c\omega) W^{\ast}.
\end{align*}

\item For $n \in \mathbb{N},$ $a,b,c \in \mathbb{C}$
\[
(D(a,b,c))^n = W \,\Diag ((a+b+c)^n, (a+bw+cw^2)^n, (a+bw^2+cw)^n)
)W^{\ast},
\]
and therefore, for $t \in \mathbb{C},$
\[
e^{tD(a,b,c)} =
W \,\Diag (e^{t(a+b+c)}, e^{t(a+bw+cw^{2})}, e^{t(a+bw^{2}+cw)}
)W^{\ast}.
\]
We note that all these matrices are in $GL(3,\mathbb{C})$
and $(e^{t D(a,b,c)})^{-1}= e^{-tD(a,b,c)} = e^{tD(-a,-b,-c)}$ for
$a,b,c,t \in \mathbb{C}.$

\item Let $a,b,c,$ $t \in \mathbb{C}.$

By (iii) $e^{tD(a,b,c)} = D (\ba(t), \bb (t), \bc(t) )$ with $\begin{pmatrix}\ba(t)\\ \bb (t)\\ \bc(t) \end{pmatrix}
= \frac{1}{\sqrt{3}} W^{\ast} \begin{pmatrix} e^{t(a+b+c)}\\
e^{t(a+b\omega+c\omega^{2})}\\ e^{t(a+b\omega^{2}+c\omega)} \end{pmatrix}.$

Therefore,
\begin{align*}
\ba(t)
&= \frac{1}{3} \left [e^{t(a+b+c)}+e^{t(a+b\omega+c\omega^{2})} + e^{t(a+b\omega^{2}+c\omega)} \right ], \\
\bb (t)&= \frac{1}{3} \left [e^{t(a+b+c)}+\omega ^2 \,e^{t(a+b\omega+c\omega^{2})} + \omega \,\, e^{t(a+b\omega^{2}+c\omega)} \right ]
\quad \mbox{and} \\
\bc(t) &=\frac{1}{3} \left [e^{t(a+b+c)}+\omega \,\, e^{t(a+b\omega+c\omega^{2})} + \omega ^2\,\,e^{t(a+b\omega^{2}+c\omega)} \right ].
\end{align*}
We set $\bd(t) = e^{td}.$ We note that $\ba, \bb, \bc, \bd$ are all
entire functions and $\ba (0) = 1 = \bd (0)$ whereas $\bb(0) = 0 =
\bc(0).$ Further, for $a,b,c,d,t$ all real, $\ba(t),$ $\bb(t),$
$\bc(t),$ $\bd(t)$ are all real.

\item Let $a,b,c,d \in \mathbb{C}.$ For $t \in \mathbb{C},$
\begin{align*}
\ba^{\prime}(t)
&= \frac{1}{3} [(a+b+c) e^{t(a+b+c)}+(a+b\omega+c\omega^2)e^{t(a+b\omega+c\omega^{2})}\\
&\quad +(a+b\omega^2+c\omega)e^{t(a+b\omega^{2}+c\omega)} ] ,\\
\bb^{\prime} (t)
&= \frac{1}{3} [ (a+b+c) e^{t(a+b+c)}+ \omega^2\, (a+b\omega+c\omega^2)e^{t(a+b\omega+c\omega^{2})} \\
&\qquad +\omega\,(a+b\omega^2+c\omega)e^{t(a+b\omega^{2}+c\omega)} ] ,\\
\bc^{\prime}(t)
&=\frac{1}{3} [(a+b+c) e^{t(a+b+c)}+\omega\,(a+b\omega+c\omega^2)e^{t(a+b\omega+c\omega^{2})} \\
&\qquad +\omega^2\,(a+b\omega^2+c\omega)e^{t(a+b\omega^{2}+c\omega)} ]\qquad \mbox{and} \\
\bd^{\prime}(t)
&= d \be^{td}.
\end{align*}

In particular, $\ba^{\prime}(0) = a,$ $\bb^{\prime}(0) = b,$
$\bc^{\prime}(0) = c,$ and $\bd^{\prime}(0) = d.$ As a consequence,
if $\ba(t_n)$ (respectively $\bb(t_n),$ $\bc(t_n)$) are all real
for a real
sequence $(t_n)$ convergent to zero then $a$ (respectively $b,c$) is
real.

Thus in view of the last line of (v) we may say that
$\ba(t_n),$ $\bb(t_n),$ $\bc(t_n)$ are all real for a real sequence $(t_n)$
convergent to zero if and only if $a,b,c$ are all real if and only
$\ba(t),$ $\bb(t),$ $\bc(t)$ are all real for all real $t.$ A
similar statement holds for the function $\bd$ as well.

\item Let $a,b,c,d \in \mathbb{C}$ and set $\rho = \rho
[a,b,c,d].$ Then by (v) above, for $t \in \mathbb{C},$ $\brho (t) =
e^{t \rho}$ coincides with $\rho [\ba(t), \bb(t), \bc(t), \bd(t)].$
We first note that in view of (iv) above, each $\brho(t)$ is a
bijective map on $M_3$ to itself. Further, for $b=0=c,$ $\ba(t) =
e^{at}$ and $\bb(t)=0=\bc(t),$ so that
$$\brho(t) = e^{ta} \mathcal{I}d_{D_{n}} \oplus e^{td}
\mathcal{I} d_{F_{n}} \,\,\mbox{for}\,\, t \in \mathbb{C}.$$

\begin{enumerate}[\rm(a)]
\item Items (i) and (vi) may be combined to give: $\brho(t_n)$ are all
$\ast$-maps for a real sequence $(t_n)$ convergent to $0$ if and only if
$a,b,c,d$ are all real if and only if $\brho(t)$ are all $\ast$-maps for
all real~$t.$

From now onwards we consider only real $a,b,c,d,t.$

\item By 2.H(ii)(c) and (v) above, for any $t \neq 0,$ $\brho(t)$ is
unital if and only if $a+b+c=0$ and in that case all $\brho(t)$ are
unital as $t$ varies in $\mathbb{R}.$ Similar statements hold with
unital replaced by trace-preserving.
\end{enumerate}

\item Let $a,b, c, d$ be real. Set $u = \frac{1}{2}
(b+c),$ $v = \frac{1}{2}(b-c).$ Then
\begin{align*}
a+b+c&= a +2u,\\
a+b\omega+c\omega^2 &= a - u + i \sqrt{3} v,\\
a+b\omega^2+c\omega &= a - u - i \sqrt{3} v.
\end{align*}
Let $t \in \mathbb{R}.$ Then
\begin{align*}
\ba(t)
&= \frac{1}{3} \left [e^{t(a+b+c)} + e^{t(a+b\omega+c\omega^2)} + e^{t(a+b\omega^2+c\omega)} \right ] \\
&= \frac{1}{3} \left [e^{t(a+2u)} + e^{t(a-u+i\sqrt{3}v)} + e^{(a-u-i\sqrt{3}v)} \right ]\\
&= \frac{1}{3} e^{t(a-u)} \left [e^{3tu} + 2 \cos (\sqrt{3} vt)\right ],\\
\bb(t)
&= \frac{1}{3} \left [e^{t(a+2u)} + e^{t(a-u+i\sqrt{3}v)-\frac{2}{3} \pi i} + e^{t(a-u-i\sqrt{3}v)+\frac{2}{3} \pi i} \right]\\
&= \frac{1}{3} e^{t(a-u)} \left [e^{3tu} + 2 \cos (\sqrt{3} vt -\frac{2}{3} \pi) \right ],\\
\bc(t)
&= \frac{1}{3} e^{t(a-u)} \left [e^{3tu} + 2 \cos (\sqrt{3}vt +\frac{2}{3} \pi) \right ],\\
\bd(t)
&= e^{td} > 0.
\end{align*}
We recall from (v) above that $\ba(0) =1,$ $\bb(0)=0=c(0),$ $\bd(0)=1.$

If $\bm{\rho}$ is a positive map, then by 2.H(ii)(c) $\ba(t),$
$\bb(t),$ $\bc(t),$ are $\ge 0.$

We begin by finding out when $\ba(t),$
$\bb(t),$ $\bc(t),$ are $\ge 0$ and then go on to find conditions
under
which $\brho(t)$ is completely positive, PPT, separable etc.

\item We can argue as in (vi) above and have that if $\bb(t_n)$
(respectively $\bc(t_n)$) are all non-negative for a sequence $(t_n)$
in $(0, \infty)$ convergent to zero, then $b$ (respectively $c$) is
$\ge 0.$ So from now onwards we take $b, c \ge 0.$

\item Let $b=0=c.$ Then $\ba(t) = e^{ta}$ and $\bd(t) = e^{td} > 0$
for all $t \in \mathbb{R}$ whereas $\bb(t) = 0 = \bc(t)$ for all $t
\in
\mathbb{R}.$ By 2.H(iii)(c) no $\brho(t)$ is completely
copositive. By 2.H(iii)(b), $\brho(t)$ is completely positive
if and only if $e^{ta} \ge e^{td}$ if and only if $ta \ge td.$
\begin{enumerate}[\rm(a)]
\item For $a = d, \{\brho (t) : t \in \mathbb{R} \} \equiv \{e^{td}
Id_{M_{3}} : t \in \mathbb{R}\}$ is a group of completely positive
maps that are all non-PPT, which illustrates Theorem~5.6. For $a = d
= 0$ it is the trivial group $\{Id_{M_{3}}\}$ for $t \in \mathbb{R}.$
\item The family $\{\brho(t) : t \ge 0\}$ is a quantum dynamical
semigroup if and only if $a \ge d$ and all the maps are non-PPT and therefore,
entangled. This illustrates the condition (iii) of the Trichotomy in
Theorem 5.2 and Remark 5.5(ii) (a).
\item Let $a = 0 > d$ and $t>0$.\ Then by (vii)(b) above each $\brho(t)$ is unital
and trace-preserving; in particular, each $\brho(t)$ is a quantum channel.
By Theorem 2.6(ii)(e), for $0\le t<- \frac{(\log 2)}{d}$. $\bm\rho(t)$ has Schmidt number 3, and is, therefore, not partially entanglement braking.

\item It follows from Remark~2.8(i)(a) that for $a=d,$ the Choi matrix
$C_{\brho(t)}$ has rank $1$ for all $t \in \mathbb{R}$ and, on the
other hand, for $a> d,$ $t>0,$ the Choi matrix $C_{\brho(t)}$ has
rank $3.$
\end{enumerate}

\item Let $(b,c) \neq (0,0),$ $b,$ $c \ge 0.$ We refer to (viii) above.

Then $u > 0,$ $u \ge |v|.$ So for $t <
0,$ $e^{3tu} < 1.$ Also $2 \cos (\sqrt{3} v t + \frac{2}{3} \pi)$
assumes
value $-1$ for some $t < 0$ and thus $\bb(t) < 0.$ Similar
conclusions hold for $\bc(t).$ So we consider only the case $t \ge
0.$ As already noted in (v) $\ba(0) =1,$ $\bb(0)=0,$ $\bc(0)=0.$

\begin{enumerate}[\rm(a)]
\item In case $b=c,$ i.e., $v=0,$ we immediately have for $t > 0,$
\begin{align*}
\ba(t) &= \frac{1}{3} e^{t(a-u)} \left [e^{3ut} +2 \right ] > 0, \\
\bb(t) &= \frac{1}{3} e^{t(a-u)} \left [e^{3ut} -1 \right ] > 0,\\
\bc(t) &= \frac{1}{3} e^{t(a-u)} \left [e^{3ut} -1 \right ]= \bb(t)
> 0.
\end{align*}
For the general case some computations are needed.
\item Let $\alpha = 0,$ $\frac{2 \pi}{3},$ $\frac{-2 \pi}{3}.$
Set $f_{\alpha}(t) = e^{3ut} + 2 \cos (\sqrt{3} vt + \alpha),$ $t \in
\mathbb{R}.$ Then $f_{\alpha}$ is infinitely differentiable and
$f_{\alpha}(0)= 1 + 2 \cos \alpha \ge 0.$ Further, for $t \in
\mathbb{R},$ $f_{\alpha}^{\prime}(t)= 3ue^{3ut} - 2 \sqrt{3} v \sin
(\sqrt{3} v t + \alpha).$ Therefore, for $t \in \mathbb{R},$
$f_{\alpha}^{\prime \prime}(t) = (3u)^2 e^{3ut} - 2
(\sqrt{3} v)^2 \cos (\sqrt{3} v t + \alpha) \ge 9 u^2 e^{3ut} - 6 v^2
= 9 u^2 (e^{3ut}-1) + (9 u^2 - 6 v^2).$ So for $t > 0,$
$f_{\alpha}^{\prime \prime}(t) > 0.$

As a consequence $f_{\alpha}^{\prime}$ is strictly increasing on $[0,
\infty).$ Now $f_{\alpha}^{\prime}(0) = 3 u - 2 \sqrt{3} v \sin
\alpha,$ which is $3u,$ $3u - 3v,$ $3u+3v$ respectively for $\alpha =
0,$ $\frac{2 \pi}{3},$ $\frac{-2 \pi}{3}$ respectively i.e. $3u,$
$3c,$ $3b$ respectively. But $3u,$ $3b,$ $3c$ are all $> 0.$ So
$f_{\alpha}^{\prime} (t) > 0$ for $t > 0.$ Therefore, $f_{\alpha}$ is
strictly increasing on $[0, \infty).$ Consequently $f_{\alpha}(t) > 0$
for
$t>0$ and hence $\ba(t),$ $\bb(t),$ $\bc(t)$ are all $> 0$ for $t>0.$

\item Now $\ba(t) \ge \bd(t)$ if and only if $\frac{1}{3} e^{t(a-u)} \left [e^{3ut}
+
2 \cos (\sqrt{3} vt)\right ] \ge e^{td}$ if and only if
$e^{-ut} \left [e^{3ut} + 2 \cos (\sqrt{3} vt)\right ] \ge 3
e^{t(d-a)}.$
Set $w = a-d.$
The condition $\ba(t) \ge \bd(t)$ is equivalent to
$$e^{-ut} \left [e^{3ut} + 2 \cos (v\sqrt{3} t) \right ] \ge
3e^{-wt}.$$
Let $g(t) = e^{2ut} + 2 e^{-ut} \cos (v\sqrt{3} t) - 3 e^{-wt}, \,\,t
\in \mathbb{R}.$ Then $g$ is infinitely differentiable and $g (0) =
0.$ Also for $t \in \mathbb{R},$
\begin{align*}
g^{\prime}(t)
&= 2 ue^{2ut} + 2 e^{-ut} \left (-u \cos (v\sqrt{3}t)- v\sqrt{3} \sin (v \sqrt{3}t)\right ) + 3 w e^{-wt}\\
&= 2 u \left (e^{2ut} - e^{-ut} \cos (v \sqrt{3}t) \right ) - 2 v\sqrt{3} e^{-ut} \sin (v \sqrt{3}t) + 3 we^{-wt} \\
&= 2u \left [(e^{2ut}- e^{-ut})+ 2e^{-ut} \sin^2(\frac{v\sqrt{3}}{2}t)\right ] - 2 v \sqrt{3} e^{-ut} \sin (v \sqrt{3}t) + 3 we^{-wt}.
\end{align*}
In particular, $g^{\prime} (0) = 3 w.$ So if $g(t_n) \ge 0$ for a
sequence $(t_n)$ in $(0, \infty)$ with $t_n$ convergent to $0$ then
$g^{\prime}(0) \ge 0,$ i.e., $w \ge 0.$ Now assume $w \ge 0.$ Then for
$t > 0,$ using $|\sin t| \leq |t|$ for all $t,$
\begin{align*}
g^{\prime}(t)
&\ge 2 u (e^{2ut} - e^{-ut}) - 2 v \sqrt{3} e^{-ut} (v\sqrt{3} t) \\
&= 2 e^{-ut} \left [u (e^{3ut} -1) - 3v^2 t \right ].
\end{align*}
Let $h(t) = u (e^{3ut}-1) - 3 v^2 t,$ $t \in \mathbb{R}.$ Then $h$ is
infinitely differentiable and $h(0) = 0.$ Also for $t> 0,$
\begin{align*}
h^{\prime}(t)
&= u. 3ue^{3ut}-3v^2 \\
&= 3 u^2 (e^{3ut}-1) + 3 (u^2-v^2) > 0\,.
\end{align*}
So $h(t) > 0$ for $t > 0.$ As a consequence, $g^{\prime}(t) >0$ for
$t> 0.$ This gives $g(t) > 0$ for all $t>0.$ Thus $g(t_n) \ge 0$ for
a sequence $(t_n)$ in $(0,\infty)$ with $t_n$ convergent to zero if and only if
$w \ge 0$ if and only if $g(t) > 0$ for all $t > 0.$ Hence $\brho(t_n)$ are all
completely positive maps for a sequence $(t_n)$ in $(0,\infty)$
convergent to $0$ if and only if $w \ge 0$ if and only if $\brho(t)$ are all completely
positive maps for all $t > 0.$

\item Moreover, Remark~2.8 then gives that for $a \ge d,$ $t>0,$ the Choi
matrix
$C_{\brho(t)}$ has rank $9.$
Consider this case together with $a+b+c=0$, i.e., $a=-2u$. We refer to Theorem 2.6((ii)(e) and Theorem 2.7(ii)(e) with $a,b,c,d$ there replaced by $\bm a(t)$, $\bm b(t)$, $\bm c(t)$, $\bm d(t)$ respectively. The inequalities $\bm d(t)\le \bm a(t) < \bm a(t)+2(\bm b(t)+\bm c(t)) <2\bm d(t)$ are satisfied for
$t$ in the non-empty interval $[0,-(\log 3-\log 2)/d]$,
for sure. So, for all such $t$, in view of (vii)(b) above each $\bm \rho(t)$ is a quantum
channel that is not  partially entanglement breaking and is non-PPT as well.

\item Suppose $w \ge 0$; then
\begin{align*}
h(t)&=\bb(t) \bc(t) - \bd(t)^2\\
&=\frac{1}{9} e^{2t(a-u)} \left [(e^{3ut} - \cos (\sqrt{3}vt))^2 - 3 \sin^2 (\sqrt{3} vt) \right ] - e^{2dt}\\
&= \frac{1}{9} e^{2t(a-u)} \left [e^{6ut} - 2e^{3ut} \cos (\sqrt{3}vt) + \cos^2 (\sqrt{3} vt) - 3 \sin^2 (\sqrt{3}vt) \right ] -e^{2dt}\\
&= \frac{1}{9} e^{2t(a-u)} \left [e^{6ut} - 2e^{3ut} \cos (\sqrt{3}vt) - 1 + 2 \cos (2\sqrt{3} vt) \right ] -e^{2dt}\\
&= \frac{1}{9} e^{2ta} \left [e^{4ut} - 2e^{ut} \cos (\sqrt{3}vt) -e^{-2ut} + 2 e^{-2ut} \cos (2\sqrt{3}vt) - 9 e^{-2wt} \right ]\\
&= \frac{1}{9} e^{2ta} \, g(t),
\end{align*}
where
$g(t) =e^{4ut} - e^{-2ut} - 9 e^{-2wt} -2e^{ut} \cos
(\sqrt{3}vt) +2e^{-2ut} \cos (2\sqrt{3}vt), t \in \mathbb{R}.
$ We note that $g$ is infinitely differentiable on $\mathbb{R}$ and
$g(0) = -9.$

Since $u > 0,$ $g(t) \rightarrow \infty$ as $t \rightarrow \infty.$
So there
is an $s_0 \in (0, \infty)$ satisfying $g(t)>0$ for $t>s_0.$
This, in turn, gives that $h(t)>0$ for $t>s_0$.
By 2.H(iii)(c) $\brho(t)$ is PPT for $t> s_0.$ As $T_0=Id$ is not
PPT, an application of the Trichotomy result as envisaged in Remark~5.5(ii)
(a) immediately gives that there exists a unique $t_0 \in (0,
\infty)$ such that for $0 \leq t < t_0,$ $\brho(t)$ is not PPT but
for $t \ge t_0,$ $\brho(t)$ is PPT. This, in view of 2.H(iii)(c),
entails that there exists a $t_0 \in (0, \infty)$ satisfying,
$h(t) < 0$ for $0 \leq t < t_0$ and $h(t) \ge 0$ for $t \ge t_0.$

We now proceed to refine this observation.
\item For $t \in \mathbb{R}$,
\begin{align*}
g^{\prime}(t)
&= 4 ue^{4ut} + 2ue^{-2ut} + 18 we^{-2wt}- 2e^{ut} \left (u \cos (\sqrt{3} vt) - \sqrt{3} v \sin (\sqrt{3}vt) \right ) \\
&\qquad + 2 e^{-2ut} \left (-2u \cos (2 \sqrt{3}vt) - 2 \sqrt{3} v \sin (2\sqrt{3} vt) \right ) \\
&= 2 u \left [2 e^{4ut} + e^{-2 ut} - e^{ut} - 2 e^{-2 ut} \right ] +18 we^{-2wt} \\
&\qquad + 4 u e^{ut} \sin^{2} (\frac{\sqrt{3}}{2} vt) + 2 \sqrt{3}ve^{ut} \sin (\sqrt{3}vt)\\
&\qquad + 8 u e^{-2ut} \sin^2 (\sqrt{3} vt) - 4 \sqrt{3}ve^{-2ut} \sin(2\sqrt{3} vt).
\end{align*}
We note that $g^{\prime}(0) =18 w \ge 0.$ Now for $t \ge 0,$
\begin{align*}
g^{\prime}(t)
&\ge 2 u \left [2e^{4ut} - e^{ut} - e^{-2ut} \right ] -2 \sqrt{3} ve^{ut} (\sqrt{3} vt)- 4 \sqrt{3} v e^{-2ut}(2\sqrt{3}vt)\\
&= 2 \left [ u (2e^{4ut}-e^{ut}-e^{-2ut})-3v^2 e^{ut} t -4 \times3v^2 e^{-2ut}t \right ]\\
&= 2 e^{-2ut} \left [u (2 e^{6ut} - e^{3ut}-1) - 3 v^2 e^{3ut} t -12 v^2 t\right ] \\
&= 2 e^{-2ut} \left [u (e^{3ut}-1) e^{3ut} + u(e^{6ut}-1)-3v^2 te^{3ut} - 12 v^2 t \right ]\\
&\ge 2 e^{-2ut} \left [3u^2 t e^{3ut} + 6u^2 t - 3v^2 t e^{3ut} - 12v^2 t \right ] \\
&= 2 e^{-2ut} \left [\frac{3}{2} u^2 t e^{3ut} + (u^2 -2v^2) \left (\frac{3}{2} t e^{3ut} + 6t \right ) \right ].
\end{align*}
Because $u>0,$ we have for $t>0,$ $g^{\prime}(t) > 0$ in case $u^2\ge 2 v^2.$

One can obtain $g^{\prime}(t) > 0$ for $t > 0$ for less restricted
cases but we prefer to confine our attention to this simple case and
go on with the case $u \ge \sqrt{2} |v|.$ Then $g$ is strictly
increasing on $[0,\infty).$ So there
exists a unique $t_0\in (0, \infty)$ such that $g(t_0) =0,$ $g(t)<0$
for $0 \leq t < t_0$ and $g(t) > 0$ for $t_0 < t < \infty.$
As a consequence, there exists a unique $t_0\in (0,\infty)$ such that
$h(t_0)=0$, $h(t)<0$ for $0\le t<t_0$ and $h(t)>0$ for $t_0<t<\infty$.
So by 2.H(iii)(c), $\brho(t)$ is not completely co-positive for $t
< t_0$ but is completely co-positive for $t \ge t_0.$

\item Hence for $(b,c) \neq (0,0),$ $a \ge d,$ $b+c \ge
\sqrt{2}|b-c|,$ there
exists a unique $t_0 \in (0,
\infty)$ that satifies
\begin{enumerate}
\item[($\alpha$)] for $0 \leq t < t_0,$ $\brho (t)$ is not PPT, and
\item[($\beta$)] for $t \ge t_0,$ $\brho (t)$ is PPT.
\end{enumerate}
\end{enumerate}

This illustrates the condition (ii) of Trichotomy in Remark~5.5(ii)(a) in a concrete manner.

\item Let $\bm{\tau}(t)=e^{t\tau[a,b,c,d]}$, $t\ge
0$. Then
$\bm{\tau}(t)=D(\bm{a}(t),\bm{b}(t),\bm{c}(t))\oplus (\cosh(td)I_{F_n}+\sinh(td)\tau_{F_n})$. Its Choi matrix in expanded form is
$$
C_{\bm{\tau}(t)}=\scriptsize
\begin{bmatrix}
\bm{a}(t) &0&0&0&\cosh(td) &0&0&0&\cosh(td)\\
0&\bm{c}(t) &0&\sinh(td) &0&0&0&0&0\\
0&0&\bm{b}(t) &0&0&0&\sinh(td) &0&0\\
0&\sinh(td)&0&\bm{b}(t) &0&0&0&0&0\\
\cosh(td)&0&0&0&\bm{a}(t) &0&0&0&\cosh(td)\\
0&0&0&0&0&\bm{c}(t) &0&\sinh(td)&0\\
0&0&\sinh(td)&0&0&0&\bm{c}(t) &0&0\\
0&0&0&0&0&\sinh(td)&0&\bm{b}(t) &0\\
\cosh(td)&0&0&0&\cosh(td)&0&0&0&\bm{a}(t)\\
\end{bmatrix}.
$$
It has trace $\mu(t)=3$ $(\bm a(t)+\bm b(t) +\bm c(t))>0$.
Further, for $t > 0,$ $C_{\bm{\tau}(t)}$ is a positive
matrix if and only if $\bm{a}(t) \ge \cosh (td),$ $\bm{b}(t) \ge 0,$
 $\bm{c}(t)\ge 0$ and $\bm{b}(t) \bm{c}(t) \ge \sinh^2 (td).$
Computations of the type done in this example give that this happens
for all $t \ge 0$ if $\frac{2}{3}b=\frac{2}{3} c \ge a \ge |d|;$ and,
in fact, for less restricted cases as well.

\item \textbf{An interesting event.}\ Around the time of acceptance of the paper, the author visited The Institute of
Quantum Computation at Waterloo University, Canada from October 25 to 30, 2015. During
a discussion, Vern Paulsen told her about the PPT  Square conjecture of Matthias Christandl which can be found as
(Problem G). In Banff International Research Station workshop: Operator structures in quantum information theory (2012).
Available at
http://www.birs.ca/workshops/ 2012/12w5084/report12w5084.pdf.
She thought her Example~5.7 would settle the conjecture in the negative.
She got in touch with Matthias Christandl and an intense discussion took place.
This led to addition of (ii)(e)  in Theorem 2.6, deletion of Schmidt number part in earlier Remark 2.9 and a consequent change in part (xi)(d) and (xii) above. It turned out that the example does not affect the conjecture so far.

\end{enumerate}
\end{example}

\medskip
{\bf Acknowledgment.}

I express my deep sense of gratitude to Kalyanapuram Rangachari
Parthasarathy. I have learnt most of the concepts in this paper from
him during his Seminar Series of Stat. Math. Unit at the Indian
Statistical Institute, New Delhi, University of Delhi and elsewhere.
I~have gained immensely from insightful discussion sessions with him
from time to time.

I thank B. V. Rajarama Bhat and Aurelian
Isar for their kind comments and suggestions. I also thank Rajendra
Bhatia for his discussion and seminar sessions which strengthened my
knowledge and interest in matrices.

It is my pleasure to thank Kenneth A. Ross for his reading of
the paper, his useful suggestions and encouragement.

Improvement, particularly in item 2.D and the general level of the
paper was motivated by quick constructive critical comments and kind
suggestions by Mary Beth Ruskai on an earlier version of the paper.
She made me aware of a good amount of relevant research work which
led to my further work. I am grateful to her for her kind advice.
In the same vein, I thank S.-H. Kye, S.R.~Garcia and S.K. Jain for
their help and encouragement in various
forms and at different stages of this final version of the article
beginning with its first version [80] (a fore-runner to this paper)
by referring to it, commenting
on it or making useful suggestions.

I thank the referees for their useful comments and suggestions that improved the paper.

In view of Example 5.7 (xiii) above, it is my pleasant duty to thank John Watrous for his kind invitation to
visit his Institute of Quantum Computing, Waterloo University, Canada and the institute
for kind hospitality. I thank him and his colleagues, particularly, Richard Cleve and Vern
Paulsen for useful discussion that improved my perspective. I thank Matthias Christandl for useful discussion, his
careful reading of the paper and suggestions for corrections and improvement in the paper.

I thank K. Manjunatha Prasad for inviting me to give a talk at the excellent
conference ICLAA 2014 and contribute an original article to the Proceedings
and keeping patience with me while processing it.

I thank Mr. Anil Shukla and M/s Scientific Documentations for their efficient typing and
cooperative
attitude in coping with my bad handwriting and odd style.

I thank Indian National Science Academy for
support under the INSA Senior Scientist and Honorary Scientist Programme and Indian
Statistical Institute, New Delhi
for  visiting positions under
this programme together with excellent research facilities.

\newpage

\vskip4mm

\section*{Appendix}

\begin{center}
\centerline{\bf Happy 60th Birthday Ravindra B. Bapat}
\begin{tabular}{l}
Spreading his special knowledge rays,\\
\mbox{}\qquad bright and colourful in a line,\\
And concentrating on inherent arrays,\\
\mbox{}\qquad creates Linear models of many a kind.\\
Looking so sober and systematic,\\
\mbox{}\qquad but goes for a random walk on a tree.\\
Makes me wonder on a change like this,\\
\mbox{}\qquad but soon he comes back with a glee.\\
Hands full of Laplacian rays,\\
\mbox{}\qquad beautiful fans with a glow,\\
Wheels moving with different speeds,\\
\mbox{}\qquad displaying wonderful life flow!\\
Felicitations to you, Ravindra Bapat,\\
\mbox{}\qquad best wishes for laurels many more,\\
Happy long life to you and your family,\\
\mbox{}\qquad and each one on ICLAA-14 floor!
\end{tabular}
\end{center}

\bigskip
\bigskip
\begin{quote}
\textbf{Abstract} of the actual expository talk at the International Conference on Linear Algebra and its Applications-2014 at Manipal University,
Mangalore, India.
\end{quote}

\begin{center}
The passage from two to three and three to infinity\\ in classic to quantum channels.
\end{center}

\bigskip
The situation changes drastically for matrices, maps on matrix algebras and
applications to Quantum Information theory when we go from
order two to three or three to infinity.
Examples include maximally entangled bases and
Quantum Birkhoff Theory.
\end{document}